\documentclass[
superscriptaddress,
amsmath,amssymb,
aps,
pre,
]{revtex4-2}

\usepackage{graphicx}
\usepackage{xcolor}
\usepackage{amsthm}
\usepackage{comment}
\usepackage{amsfonts}
\usepackage{nccmath}

\providecommand{\U}[1]{\protect\rule{.1in}{.1in}}

\newtheorem{theorem}{Theorem}

\newtheorem{lemma}[theorem]{Lemma}

\newtheorem{remark}[theorem]{Remark}

\newcommand{\R}{\mathbb{R}}

\newcommand{\mfe}{\mathfrak{e}}

\newcommand{\EE}{\mathbb{E}}

\newcommand{\m}[1]{\begin{pmatrix}#1\end{pmatrix}} 
\newcommand{\avg}[1]{\langle#1\rangle} 
 
\newcommand{\dd}{\textup{d}}

\begin{document}

\preprint{preprint}

\title{Sensitivity analysis of colored noise-driven interacting particle systems}

\author{Josselin Garnier}
 \email{josselin.garnier@polytechnique.edu}
 \affiliation{Centre de Mathématiques Appliquées, Ecole Polytechnique, Institut Polytechnique de Paris, 91120 Palaiseau, France \\}
\author{Harry L. F. Ip}
 \email{harryip2@cityu.edu.hk}
\author{Laurent Mertz}
 \email{lmertz@cityu.edu.hk}
 \affiliation{Department of Mathematics, City University of Hong Kong, \\ Kowloon, Hong Kong, China \\}

\date{\today}

\begin{abstract}
We propose an efficient sensitivity analysis method for a wide class of colored noise-driven interacting particle systems (IPS). 
Our method is based on unperturbed simulations and significantly extends the Malliavin weight sampling method proposed by Szamel (EPL, 117 (2017) 50010) for evaluating sensitivities such as linear response functions of IPS driven by simple Ornstein-Uhlenbeck processes. We show that the sensitivity index depends not only on two effective parameters that characterize the variance and correlation time of the noise, but also on the noise spectrum. In the case of a single particle in a harmonic potential, we obtain exact analytical formulas for two types of linear response functions. By applying our method to a system of many particles interacting via a repulsive screened Coulomb potential, we compute the mobility and effective temperature of the system. Our results show that the system dynamics depend, in a nontrivial way, on the noise spectrum. 
\end{abstract}

\maketitle

\tableofcontents

\section{Introduction}
Sensitivity analysis methods were initially developed for stochastic processes used in financial engineering, specifically in dealing with price sensitivities of financial products - often referred to as “Greeks” (see, for instance, \cite[Chapter 7]{glasserman2004monte}
and \cite{chen2007malliavin}). 
Over the last decade, applications of these methods have emerged in the physics literature for non-equilibrium systems. They address the sensitivities of complex, white, and colored noise-driven systems (e.g., active matter) with respect to certain types of perturbations
\cite{Szamel_2017, 
warren2012malliavin,dal2019linear}. In this area, a fundamental problem is to characterize the sensitivity of complex systems to perturbations by an external force from their steady state. Understanding the changes of a physical system in response to an internal or external perturbation can give crucial insights into the underlying physics, as evidenced by \cite{dal2019linear,  das2023nonequilibrium, fodor2016far, wittmann2017effective, PhysRevE.100.022601, PhysRevE.107.034110}.
Our goal is to show that it is possible to answer these questions from unperturbed simulations even when the system and the driving noise are complex.
It should also be pointed out that sensitivity indices, such as linear response functions for equilibrium systems, can be obtained
from the fluctuation-dissipation relations (FDRs) and correlation functions evolving with unperturbed dynamics \cite{chandler1987introduction,hansen2013theory}. We further mention that techniques that allow one to calculate the sensitivity indices with simulations of the unperturbed system have been proposed \cite{chatelain2003far,ricci2003measuring,corberi2010fluctuation} in a different context, that of Ising spin systems (discrete in space or time).

\subsection{Review of fundamental work}
A considerable leap forward in this direction of research has been made with the work of Warren and Allen \cite{warren2012malliavin}. They introduce a simple computational approach for responses to infinitesimal changes in internal or external parameters in stochastic Brownian dynamics simulations, which may be in or out of equilibrium.
The method is closely related to the methods developed for the ``Greeks'', particularly those related to likelihood ratios and Malliavin calculus.
Indeed, the method does not require simulating the perturbed system, but in addition to simulating the unperturbed system, it requires following extra stochastic variables. The latter corresponds to the derivatives of the probability density with respect to the parameter of interest.  Warren and Allen have coined the term ``Malliavin weight'' for these variables. They consider a cluster of $N$ particles, each in $\mathbb{R}^d$, $x_\lambda = (x_{1,\lambda}, \ldots, x_{N,\lambda})$ satisfying 
\begin{equation}
\label{eq:dedowa}
\dot x_\lambda =  \frac{D}{k_B T_{\textup{eff}}^{\textup{sp}}}  F_\lambda(x_\lambda) +   \sqrt{2D} \dot w,
\end{equation}
where $\dot w$ is an $n$-dimensional Gaussian white noise (time derivative of a $n$-dimensional Wiener process $w$),
$n = Nd$, 
$D>0$ is the diffusion coefficient (assumed to be the same for all particles), $k_B$ is the Boltzmann constant, and $T_{\textup{eff}}^{\textup{sp}}$ is the single-particle effective temperature. Here, $\lambda$ is a scalar (real) parameter of interest for sensitivity analysis and the function $F_\lambda = (F_{1,\lambda}, \ldots, F_{N,\lambda})$ specifies the external forces.
The Malliavin weight is the real-valued random process $q_\lambda$ given by 
\begin{equation}
\label{eq:mwwa}
\forall t \geq 0, 
\: \: q_\lambda(t) = \frac{\sqrt{D}}{\sqrt{2} k_B T_{\textup{eff}}^{\textup{sp}}} \int_0^t \frac{\partial F_\lambda}{\partial \lambda}(x_\lambda(s)) \cdot \dd w (s).
\end{equation}
The notation $\cdot$ is the inner product in $\mathbb{R}^{n}$, and the equation should be understood as an It\^o's stochastic integral.   
Warren and Allen have shown that for any test function $\Phi:\mathbb{R}^n \to \mathbb{R}$
\begin{equation}
\label{FORMULA1}
\forall t \geq 0, 
\: \: \frac{\dd}{\dd \lambda} \langle \Phi(x_\lambda(t)) \rangle_\lambda  \bigg|_{\lambda = 0}
=
\langle  \Phi(x_0(t)) q_{0}(t) \rangle.
\end{equation}    
They have applied their method to a nonequilibrium-driven steady state formed by a cluster of particles in a two-dimensional harmonic trap under shear.  Here, $\langle \cdots \rangle_\lambda$ denotes averaging for the system prepared at $t=0$ in the steady state corresponding to the force $F_0$ and then evolving for $t>0$ under the influence of the modified force $F_\lambda$. The notation $\langle \cdots \rangle$ is for the unperturbed steady state average $\langle \cdots \rangle_0$. 
Eq.~(\ref{FORMULA1}) shows that the estimation of the sensitivity index can be carried out with trajectories of the unperturbed dynamics (\ref{eq:dedowa}) at $\lambda=0$,  only, together with the Malliavin weight (\ref{eq:mwwa}) at $\lambda=0$.
Warren and Allen's technique has been employed in \cite{dal2019linear} to test analytical predictions with new Green-Kubo-like expressions for the diffusivity and mobility with unperturbed simulations. It has also been used in \cite{petrelli2020effective} to complete an experimental study on the effective temperature via the violation of FDRs for a system of active Brownian particles (ABP) with inertia. The authors have investigated theoretically and numerically the linear response and effective temperature for a single particle.
It should be noted that beyond the scope of sensitivity analysis,
significant applications of Formula \eqref{FORMULA1} have been proposed for optimizing the steady state of a self-assembling colloid \cite{das2023nonequilibrium}.
This formula makes it possible to obtain an explicit formula for gradients used in the optimization procedure, thus avoiding errors in approximations such as with finite differences \cite{das2019variational}.

Another success in this line of research is a nontrivial generalization of the Malliavin Weight Sampling (MWS) method as proposed by Szamel \cite{Szamel_2017}.
He introduced a method for calculating sensitivities of statistics for non-equilibrium systems expressed as colored noise-driven systems. In particular, he studied the evolution of  a cluster of particles under the influence of self-propulsion represented by a colored noise. 
The components of the particles are driven by independent, identically distributed stationary Gaussian noises with a correlation time $\tau_p > 0$. The driving process $f$ (taking values in $\mathbb{R}^n$, $n=Nd$) is of the form
\begin{equation}
\label{simple_ou}
\tau_p \dd f = -f \dd t +   \sigma  \dd w ,
\end{equation}
where $w$ is a $n$-dimensional Wiener process  and $\sigma>0$ is the noise strength. More precisely, $\sigma = \sqrt{2 \xi_0 T_{\textup{eff}}^{\textup{sp}}}
= \xi_0 \sqrt{2 D_0}
$
where $\xi_0$ is the viscous friction coefficient and $D_0 = T_{\textup{eff}}^{\textup{sp}} \xi_0^{-1}$.
The resulting state variable describing  the cluster of particles $x$ (taking values in $\mathbb{R}^n$) satisfies 
\begin{equation}
\label{eq:ips}
\dot x =   \xi_0^{-1}(F(x) + f).
\end{equation}
Since Szamel is interested in sensitivities  (under steady state) in response to drift perturbations, he replaces $F$ by $F_\lambda = F + \lambda \hat{F}$ for some function $\hat{F}$ and he studies the derivative with respect to $\lambda$ at $\lambda = 0$ of some moments. 
With $F_\lambda$, the state variable is denoted by $x_\lambda^{\tau_p}$. The sensitivity of its statistics can be computed using the following formula valid for any test function $\Phi:\mathbb{R}^n  \to \mathbb{R}$:
\begin{equation}
\label{FORMULA2}
\forall t \geq 0, 
\: \: \left. \frac{\dd}{\dd \lambda} \langle \Phi(x_\lambda^{\tau_p}(t)) \rangle_\lambda \right |_{\lambda = 0} 
=
\langle \Phi(x_0^{\tau_p}(t)) (q_0^{\tau_p}(t) + p_0^{\tau_p}(t)) \rangle  
+ \tau_p \langle \dot \Phi(x_0^{\tau_p}(t)) q_0^{\tau_p}(t) ) \rangle  ,
\end{equation}
where the auxiliary variables $p_0^{\tau_p}$ and $q_0^{\tau_p}$ are called Malliavin weights and satisfy 
\begin{align*}
q_0^{\tau_p}(t) & = 
\frac{1}{\xi_0 \sigma} \int_0^t \hat{F} (   x_0^{\tau_p}(s) ) \cdot \dd w(s), \\
p_0^{\tau_p}(t) & = 
\frac{\tau_p}{\xi_0^2 \sigma} \int_0^t \Big[ \Big( \big(F(x_0^{\tau_p}(s))+f(s)\big) 
\cdot  \nabla_x \Big) \hat{F}
(x_0^{\tau_p}(s)) \Big] \cdot \dd w(s).
\end{align*}
Here, 
$\dot \Phi(x_0^{\tau_p}(t))$ means $\dd \left ( \Phi(x_0^{\tau_p}(t)) \right ) / \dd t$. 
When $\tau_p$ vanishes in Szamel's formula (\ref{FORMULA2}), the second term and $p_0^{\tau_p}$ vanish as well, leaving a term that is equivalent to Warren and Allen's formula (\ref{FORMULA1}). 
In addition to several examples given by Szamel that we shall examine in this paper, noticeable applications of his method have been given by Maggi et al. \cite{maggi2022critical} in which they focused on the dynamics of active particles self-propelled by \eqref{simple_ou}. This system falls within the framework of Szamel, where for each $1 \leq i \leq N$, $F_i(x) = -\sum_j \nabla_{x_i} \phi(\| x_i - x_j \|)$ and $\|x_i - x_j \|$ is the distance between two particles $x_i$ and $x_j$. The repulsive interaction potential between two particles is $\phi(r) = \beta r^{-\alpha}$, $\alpha, \beta >0$.
All particles are in $\R^d$ with $d = 2$ or $3$. Such an active matter system driven out of thermal equilibrium is a fundamental model for the physics of self-propelled particles (natural or artificial) and has been studied by many authors \cite{fodor2016far,wittmann2017effective,PhysRevE.100.022601,PhysRevE.107.034110}.

\subsection{Our contribution: sensitivity analysis with general spectra}
In Section \ref{sec:formulas}, we propose new formulas for the sensitivity indices of complex systems: complex systems such as interacting particle systems modeled by an equation of the form \eqref{eq:ips}, where the driving noise (self-propulsion) $f$ is  a stationary Gaussian noise with mean zero and with a general spectrum.
We distinguish two types of noise: noise with the same regularity as Brownian motion and noise with more regularity.
In addition, we obtain, in each of these cases, formulas generalizing \eqref{FORMULA2}. 
Naturally, in the same spirit as the works of Warren \& Allen and Szamel, these formulas allow for the evaluation of sensitivities such as time-dependent linear response functions for active particle systems propelled by persistent (colored) noise from unperturbed simulations. To compare results and to clarify the role of the noise spectrum, we calibrate the parameters of the noise spectra to keep the same amplitude (temperature) and correlation time for the different types of noise. In doing so, we will be able to ascertain whether these parameters are sufficient to determine the sensitivity indices.
Furthermore, our formulas agree with Szamel's formulas when the noise takes the form of a simple Ornstein-Uhlenbeck (OU) process as in \eqref{simple_ou}. We recover Szamel's formulas from ours when the noise is a simple Ornstein-Uhlenbeck (OU) process, as in \eqref{simple_ou}. In Section \ref{sec:harmonic}, we consider a single particle in a harmonic potential and obtain exact analytical formula for two types of perturbations (linear response functions). Finally, in Section \ref{sec:ips}, we apply our method to compute the mobility function and effective temperature of a system of many particles interacting via a repulsive screened Coulomb potential.
Our study shows that the sensitivity indices also depend on the noise spectrum, not only on the two effective parameters that characterize the amplitude and correlation time of the noise.
Coupled with the general need to accurately capture the physical impacts of perturbations on interacting particle systems driven by colored noise, these findings demonstrate the usefulness and relevance of our sensitivity analysis method.

\section{Noise structures}
\label{sec:noise}
In Subsection \ref{subsec:specgen}, we give the general form of driving noise for which we propose new representation formulas for the sensitivity indices in Section \ref{sec:formulas}. 
In Subsection \ref{subsec:specpart}, we provide a few explicit examples of noise spectra that turn out to give different sensitivity indices despite having identical variance and correlation time.

\subsection{General formulation}
\label{subsec:specgen}
We consider the system for the state variable $x$ driven by the noise $f=Cy$:
\begin{align}
\label{eq:x}
\dot x &= \xi_0^{-1} \left ( F(x) + C y \right ), \\
\dot y &= - A y + B \dot w.
\label{eq:y}
\end{align}
This is the unperturbed dynamics. 
The perturbed dynamics are obtained by replacing $F$ with $F_\lambda$.
Here 
\begin{itemize}
\item
The state variable $x(.)$ takes values in $\mathbb{R}^n$, the driving noise process $y(.)$ takes values in $\mathbb{R}^q$, and the Wiener process $w(.)$ is $p$-dimensional, $n,q,p \geq 1$,
\item
$F_\lambda = F + \lambda \hat F$ for $\lambda \in \mathbb{R}$, where $F,\hat{F}$ are functions from $ \mathbb{R}^n$ to $\mathbb{R}^n$, 
\item 
$A,B,C$ are matrices:
$A \in \mathbb{R}^{q \times q}$,
$B \in \mathbb{R}^{q \times p}$, 
$C \in \mathbb{R}^{n \times q}$.
\end{itemize}
The matrices $A$, $B$, $C$ and the function $\hat{F}$ are assumed to satisfy the following hypotheses:
\begin{itemize}
\item
All eigenvalues of $A$ have positive real parts (equivalently, $-A$ is a 
stable matrix).
\item
The matrix $B B^T$ is nonsingular, or $BB^T$ is singular but $A$, $B$, and $C$ have the Brunowski form \cite{brunovsky1970classification}:
There exist integers $n', q' \geq 1$, matrices $A_k  \in \mathbb{R}^{q'\times q'}$ for $k=1,\ldots,n'$, a matrix $\bar{B} \in \mathbb{R}^{q'\times p}$, and a matrix $\bar{C} \in \mathbb{R}^{n \times q'}$,
such that  $q = n' q'$,
\begin{align} \label{eq:ABthem2}
    A & = \m{
        O_{q'} & -I_{q'} & O_{q'} & \dots & O_{q'} \\
        O_{q'} & O_{q'} & -I_{q'} & \cdots & O_{q'} \\
        \vdots & \vdots & \vdots & \ddots & \vdots \\
        O_{q'} & O_{q'} & O_{q'} & \cdots & -I_{q'} \\
        A_1 & A_2 & A_3 & \cdots & A_{n'}
    } , 
    \quad
     B = \m{
        O_{q'\times p} \\
        O_{q' \times p} \\
        \vdots \\
        O_{q' \times p} \\
        \bar B
    } , 
    \quad \nonumber \\
    C & = \m{ \bar{C} & O_{n\times q'} & \cdots & O_{n \times q'}},
\end{align}
and $\bar B \bar B^T \in \mathbb{R}^{q' \times q'}$ is nonsingular (here $O_{p\times q}$ denotes the null matrix of size $p\times q$,
$O_p=O_{p\times p}$, and $I_p$ is the identity matrix of size $p$).
\item
$\hat F \in \textup{Im}(C)$, more exactly, there exists a function $ E$ from $\mathbb{R}^n$ to $\mathbb{R}^q$ 
such that $\hat F = C E$ when $BB^T $ is nonsingular, or there exists a function $ \bar E$ from $\mathbb{R}^n$ to $\mathbb{R}^{q'}$ 
such that $\hat F = \bar C \bar E$ when $A,B$, $C$ have the Brunowski form (\ref{eq:ABthem2}) and $\bar{B}\bar{B}^T$ is nonsingular.
\end{itemize}

The first condition ensures the ergodicity of the driving noise $y$. The second condition ensures that the stationary distribution has a multivariate Gaussian density  (with the asymptotic variance matrix $\int_0^\infty \exp(-As) B B^T \exp(-A^T s) \dd s$), and it is also used to get a simple enough representation formula of the sensitivity index (a more general hypothesis such as the Kalman rank condition that ``the augmented matrix $( B \: AB \: A^2B \: \cdots \: A^{q-1} B )$  has rank $q$" could be considered but it would lead to complicated formulas that go beyond the scope of this paper). 
The third condition guarantees that the derivative of $F_\lambda$ with respect to $\lambda$ belongs to the space explored by the noise. This is an important property that allows us to get an expression for a sensitivity index with respect to $\lambda$ in terms of an expectation that involves only the unperturbed trajectory at $\lambda=0$ and a Malliavin weight.
The model above is versatile, as it covers the case of one one-dimensional particle and the case of many multi-dimensional particles, as we will see below. 

In the next subsection, we present examples of colored noise with general spectrum that satisfy the above assumptions: the Gaussian processes with $n$ independent and identically distributed (i.i.d.) components with power spectral density (\ref{eq:psd1}) belong to the first case ($BB^T$ is nonsingular) while the Gaussian processes with $n$ i.i.d. components with power spectral density (\ref{eq:psd2}) belong to the second case ($BB^T$ is singular but $A,B,C$ have the Brunowski form).

\subsection{Particular colored noise structures and spectra}
\label{subsec:specpart}
For a $n$-dimensional OU process  $f$ defined by \eqref{simple_ou},  
the correlation function of each component is 
$$
c_{\rm ou}(t) = \frac{\xi_0 T_{\textup{eff}}^{\textup{sp}}}{\tau_p} \exp \left ( - \frac{|t|}{\tau_p} \right )
$$
and the corresponding power spectral density (PSD) is
\begin{equation}
\label{eq:psd0}
\hat c_{\rm ou}(\omega) \triangleq \int_{-\infty}^\infty e^{-i\omega t} c_{\rm ou}(t) \dd t = \frac{2 \xi_0 T_{\textup{eff}}^{\textup{sp}}}{1+\omega^2 \tau_p^2}.
\end{equation}
The OU process $f$ defined by (\ref{simple_ou}) belongs to the first case described in Subsection \ref{subsec:specgen} ($BB^T$ is nonsingular) because (\ref{simple_ou}) reads as (\ref{eq:y}) with $A= \tau_p^{-1} I_n$, $B=\sigma \tau_p^{-1} I_n$, and $C=I_n$.

We want to consider more general situations in which
$f$ is an $n$-dimensional process whose components are i.i.d. and belong to a class of zero-mean stationary Gaussian processes (SGP) with a PSD which we denote $\hat c(\omega)$. 
We want to compare situations with the same effective temperature  $T_{\textup{eff}}^{\textup{sp}}$, which is given in terms of the PSD by 
$$
2 \xi_0 T_{\textup{eff}}^{\textup{sp}} =\hat c(0) = \int_{-\infty}^\infty c(t) \dd t.
$$ 
For white noise forces $f = f_{\textup{wn}}$, a comparable situation is $f_{\textup{wn}} = \sqrt{2 \xi_0 T_{\textup{eff}}^{\textup{sp}}} \dot{w}$ (where $\dot{w}$ is a $n$-dimensional white noise) so that $\hat{c}_{\textup{wn}}(\omega)=2 \xi_0 T_{\textup{eff}}^{\textup{sp}}$, but here the correlation time is $0$. When dealing with persistent/colored noise with PSD $\hat{c}(\omega)$, in addition to having the same $T_{\textup{eff}}^{\textup{sp}}$, we want to compare situations with the same correlation time $\tau_c$. We propose to define it as the root-mean-square (rms) width of the correlation function through the relation 
$$
\tau_c^2 \triangleq - \frac{\hat c''(0)}{\hat c(0)} = \frac{\int_{-\infty}^\infty t^2 c(t) \dd t}{\int_{-\infty}^\infty c(t) \dd t}.
$$
For an $n$-dimensional OU process  $f$ defined by \eqref{simple_ou}, this gives $\tau_c = \sqrt{2} \tau_p$. 
We want to check whether the form of the PSD is important or whether the knowledge of $T_{\rm eff}^{\textup{sp}}$ and $\tau_c$ is sufficient to characterize the dynamics and the sensitivity of a system driven by such a colored noise. 
Although the results on the sensitivity indices will be established for general noise models,  
in the numerical applications, we will consider two special forms of PSD that we describe now.

\begin{enumerate}
\item 
First, we will consider a PSD of the form
 \begin{equation}
\label{eq:psd1}
\hat c(\omega) \triangleq \frac{1}{2}\sum_{k=1}^r \frac{\sigma_k^2}{1+(\omega-\omega_k)^2/\ell_k^2} + \frac{\sigma_k^2}{1+(\omega+\omega_k)^2/\ell_k^2}.
\end{equation}
Its corresponding correlation function is 
$$
c(t) = \frac{1}{2} \sum \limits_{k=1}^r
\sigma_k^2 \ell_k e^{-|t|\ell_k}
\cos(\omega_k t).
$$

When $r=1$,  we can consider parameters $(\sigma_1,\omega_1,\ell_1)$ satisfying
\begin{align*}
\hat c(0) = \frac{\sigma_1^2}{1+\gamma_1^2}
= 2 \xi_0 T_{\textup{eff}}^{\textup{sp}} 
\: \mbox{ and } \: 
- \frac{\hat c''(0)}{\hat c(0) }
= \frac{p_1 \sigma_1^2/\ell_1^2}{2 \xi_0 T_{\textup{eff}}^{\textup{sp}}}
= \tau_c^2,
\end{align*}
where $p_1 \triangleq 2(1-3\gamma_1^2)/(1+\gamma_1^2)^3$ and $\gamma_1 \triangleq \omega_1/\ell_1$,
so that the effective temperature is $T_{\rm eff}^{\textup{sp}}$ and the correlation time is $\tau_c$.
As an example, we will take 
\begin{equation}
\label{eq:parameters_psd1}
\sigma_1 \triangleq \frac{\sqrt{5}}{2} \sqrt{2\xi_0 T_{\textup{eff}}^{\textup{sp}}}, \qquad \ell_1 \triangleq \frac{2\sqrt{2}}{5}\tau_c^{-1}, \qquad \omega_1 \triangleq \ell_1/2,
\end{equation}
 and we then have
\begin{equation}\label{eq:psd1:part}
\hat{c}(\omega) = \frac{2 \xi_0 T_{\textup{eff}}^{\textup{sp}}  (1+5\omega^2 \tau_c^2/2)}{1+3\omega^2 \tau_c^2 + 25 \omega^2 \tau_c^4 /4}.
\end{equation}

\item
Second, for $\nu>0$, we will consider a PSD of the Mat\'ern form 
\begin{equation}
\label{eq:psd2}
\hat{c}(\omega) \triangleq \sigma_\nu^2
\frac{2 \pi^{\frac{1}{2}} \Gamma(\nu+\frac{1}{2}) (2\nu)^\nu}{\Gamma(\nu) \tau_\nu^{2\nu}}
\Big( \frac{2\nu}{ \tau_\nu^2}  + \omega^2\Big)^{-\nu-\frac{1}{2}}, 
\end{equation}
where we recall the definition of the Gamma function 
$
\Gamma(\mathfrak{z}) \triangleq \int_0^{\infty} t^{\mathfrak{z}-1} e^{-t} \dd t$.
Note that with $\nu=1/2$ we recover the PSD of an OU process similar to \eqref{eq:psd0}. 
When $\nu$ is larger, we deal with a smoother noise, whose trajectories can be differentiable up to the order $\lfloor\nu \rfloor$. We have $\hat{c}(0)=\sigma_\nu^2 \tau_\nu \sqrt{2\pi} \Gamma(\nu+1/2)/[\sqrt{\nu}\Gamma(\nu)]$, $\tau_c^2=\tau_\nu^2 (2\nu+1)/(2\nu)$.
Therefore, we have to take 
$$
\tau_\nu \triangleq \tau_c\sqrt{2\nu/(2\nu+1)} 
\: \mbox{ and } \: \sigma_\nu^2 \triangleq \xi_0T_{\textup{eff}}^{\textup{sp}} \sqrt{2\nu+1}\Gamma(\nu) /[ \Gamma(\nu+1/2) \sqrt{\pi}\tau_c],
$$
so that the effective temperature is $T_{\rm eff}^{\textup{sp}}$ and the correlation time is $\tau_c$.
As an example, we will take $\nu = 3/2$,
$\sigma_\nu^2 = \xi_0 T_{\textup{eff}}^{\textup{sp}} \tau_c^{-1}$,
$\tau_\nu = \tau_c \sqrt{3}/2$, and we then have
\begin{equation}
\label{eq:psd2:part}
\hat{c}(\omega) = \frac{2 \xi_0 T_{\textup{eff}}^{\textup{sp}}}{(1+\omega^2 \tau_c^2/4)^2}.
\end{equation}
\end{enumerate}

\begin{figure*}[t!]
\centering
\includegraphics[scale=2.15]{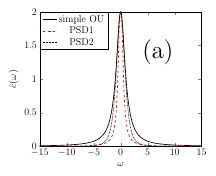}
\includegraphics[scale=2.15]{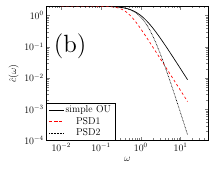}
\includegraphics[scale=2.15]{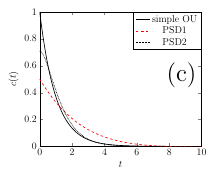}
\includegraphics[scale=2.15]{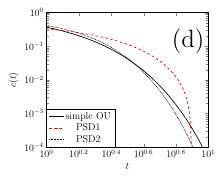}
\caption{The solid, dashed, and dotted lines represent respectively the PSD \eqref{eq:psd0} of a OU process, PSD1 \eqref{eq:psd1:part} and PSD2 \eqref{eq:psd2:part} in all the figures. Figures (a) and (b) represent the PSDs for $\omega \in [-15, 15]$ and for $\omega \in [0, 50]$ in loglog scale, respectively. Figures (c) and (d) represent the corresponding correlation functions for $t \in [0,1]$ and for $t \in [1,10]$ in loglog scale, respectively. Here $\xi_0=T_{\rm eff}^{\textup{sp}}= 1$ and $\tau_c = \sqrt{2}$. 
}
\label{fig:spectra_correlation}
\end{figure*}

In Figure \ref{fig:spectra_correlation}, we plot the PSDs \eqref{eq:psd0}, \eqref{eq:psd1:part}, and \eqref{eq:psd2:part} and the corresponding correlation functions. On the one hand, PSDs \eqref{eq:psd1:part} and \eqref{eq:psd2:part} decrease faster to $0$ than PSD \eqref{eq:psd0}, as $|\omega| \to \infty$, with PSD \eqref{eq:psd2:part} decreasing even faster than PSD \eqref{eq:psd1:part}.
On the other hand, the correlation function corresponding to \eqref{eq:psd0} decreases faster than its counterparts for \eqref{eq:psd1:part} and \eqref{eq:psd2:part} as $t \to \infty$. 

The $n$-dimensional zero-mean SGPs with power spectral density 
\eqref{eq:psd1} and \eqref{eq:psd2}
have the same distribution as the random processes $f \triangleq C y$ with $y$ taking values in $\mathbb{R}^q$
and is the solution of (\ref{eq:y}),
with $A \in \mathbb{R}^{q \times q}$,  $B \in \mathbb{R}^{q\times p}$, $C \in \mathbb{R}^{n\times q}$, and $w$ is a $p$-dimensional Wiener process. In the particular cases of the PSDs defined by (\ref{eq:psd1}) or (\ref{eq:psd2}), the matrices $A,B,C$ have the following forms.
\begin{enumerate}
\item
\begin{align}
\label{eq:matrices_for_PSD1}
A & = \oplus_{k=1}^r \m{ \ell_k I_n & -\omega_k I_n \\  \omega_k I_n &\ell_k I_n  } ,
\quad
B =  \oplus_{k=1}^r \m{  \ell_k I_n & O_n \\  O_n &  \ell_k I_n  },
\quad \nonumber\\
C & = \m{ \sigma_1 I_n & O_n & \cdots & \sigma_r I_n & O_n }.
\end{align}
\item
The PSD (\ref{eq:psd2}) with $\nu=r-1/2$ for a positive integer $r$ corresponds to $p=n$, $q=rn$,
\begin{align}
\label{eq:matrices_for_PSD2}
A & = \ell \m{
O_n  & -I_n     & O_n      &  \hdots     &   O_n  \\
O_n  & O_n & -I_n &   \hdots  & O_n \\
\vdots   & \ddots & \ddots  & \ddots &   \vdots\\
O_n &  O_n & O_n & \hdots  & -I_n \\
\alpha_0 I_n & \alpha_1 I_n   &    \cdots       & \cdots      & \alpha_{r-1} I_n
}, 
%
%
\quad 
B =\ell  \m{
O_n\\
O_n\\
\vdots \\
O_n \\
I_n
},
%
%
\quad \nonumber \\
C & = \sigma  \m{
I_n &
O_n &
\cdots &
O_n
},
\end{align}
where $\alpha_k \triangleq \binom{r}{k} = \frac{r!}{k!\,(r-k)!}$,
\begin{equation}
\label{eq:parameters_psd2}
\ell= 2 \tau_\nu^{-1} \sqrt{(2\nu+1)/(2\nu)}= 2 \tau_c^{-1}
\: \mbox{ and } \: \sigma^2 = 2 \xi_0 T_{\textup{eff}}^{\textup{sp}}. 
\end{equation}
\end{enumerate}



\section{Sensitivity indices of colored noise driven systems}
\label{sec:formulas}
We recall that $\langle \cdots \rangle_\lambda$ denotes averaging for the system prepared at $t=0$ in the steady state corresponding to the force $F(x)$ and then evolving for $t>0$ under the influence of the modified force $F_\lambda(x)$.

\subsection{The case where $B B^T$ is nonsingular}
If $(x,y)$ satisfies (\ref{eq:x}-\ref{eq:y}) with a nonsingular matrix $BB^T$, and $\Phi(.)$ is a continuously differentiable function, then
\begin{equation}
\label{eq:sensitivity1}
 \frac{\dd}{\dd\lambda} \langle \Phi(x(t)) \rangle_\lambda \Big|_{\lambda = 0}  = \Big\langle \Phi(x(t)) \big[ p_{0,0}(t) + p_{1,0}(t) \big] \Big\rangle
+ \Big\langle \dot \Phi(x(t)) \; p_{1,1}(t) \Big\rangle  ,
\end{equation} 
where  the $p_{j,k}(t)$ evolve according to the following equations of motion
\begin{equation}
\label{eq:weight1}
\begin{cases}
    & \dot p_{0,0} = E(x) \cdot A^T (BB^T)^{-1} B \dot w,
    \\& \dot p_{1,0} = \nabla E(x) \xi_0^{-1} (F(x)+C y) \cdot (BB^T)^{-1} B \dot w,
    \\& \dot p_{1,1} = E(x) \cdot (BB^T)^{-1} B \dot w,
\end{cases}
\end{equation}
with initial conditions $p_{j,k}(0)=0$ for each $0 \leq k \leq j \leq 1$ 
(here $\nabla E (x) \in \mathbb{R}^{q\times n}$ for all $x\in \mathbb{R}^n$).
Here, the notation $\cdot$ is the inner product in $\mathbb{R}^q$.
The proof is shown in Appendix \ref{appendix:proof_formula_1}. This is the first generalization of Szamel's formula.

\subsection{The case where $B B^T$ is singular but $A,B,C$ have the Brunowski form} 
If $(x,y)$ satisfies (\ref{eq:x}-\ref{eq:y}) where $A,B,C$ have the Brunowski form (\ref{eq:ABthem2}),
and $\Phi(.)$ is a function continuously differentiable at the order $n'$, then
\begin{equation}
    \label{eq:sensitivity2}
    \frac{\dd}{\dd \lambda} \avg{ \Phi(x(t)) }_\lambda \Big|_{\lambda = 0} = 
    \sum_{k=0}^{n'} \left \langle \frac{\dd^k \Phi(x(t))}{\dd t^k} \bigg[ \sum_{j=k}^{n'} \binom{j}{k} p_{j,k} (t) \bigg] \right \rangle.
\end{equation}

For $j,k \in \mathbb{N} $ such that $0\leq k\leq j\leq n'$, the function $p_{j,k}(t)$ evolves according to the following equation of motion
\begin{equation}
\label{eq:weight4formula2}
    \dot p_{j,k}(t) = \frac{\dd^{(j-k)} \bar E(x(t))}{\dd t^{(j-k)}} \cdot (A_{j+1})^T (\bar B \bar B^T)^{-1} \bar B \dot w(t),
\end{equation}
with initial condition $p_{j,k}(0)=0$ and with $A_{n'+1}:=I_{q'}$. Here, the notation $\cdot$ is the inner product in $\mathbb{R}^{q'}$, and $\binom{j}{k}$ is for $\frac{ j!}{k! (j-k)!}$. 
The proof is given in Appendix \ref{app:proof2}.
This is the second generalization of Szamel’s formula.

\begin{remark}
    Suppose we have matrices of the form $\ell A$ and $\ell B$, then the functions $p_{j,k}(t)$ become
\begin{equation}
\label{eq:secondgeneralization}
    \dot p_{j,k}(t) = \frac{\dd^{(j-k)} \bar E(x(t))}{\dd t^{(j-k)}} \cdot \frac{1}{\ell^j} (A_{j+1})^T (\bar B \bar B^T)^{-1} \bar B \dot w(t), \qquad \text{for $0\leq j \leq n'$}.
\end{equation}
\end{remark}

\begin{remark}
The first formula (\ref{eq:sensitivity1}) covers a noise structure for $C y$ that is appropriate to describe the PSD \eqref{eq:psd0} or \eqref{eq:psd1}, for instance. In general, when $B B^T$ is nonsingular, then the sample paths of $Cy$ have the same regularity as the Brownian paths. 
The second formula (\ref{eq:sensitivity2}) covers a noise structure for $Cy$ that is appropriate to describe the PSD \eqref{eq:psd2}, for instance. In general, when $A,B,C$ have the Brunowski form and $\bar{B} \bar{B}^T$ is nonsingular, then the sample paths  of $Cy$ are smooth. Indeed, this is equivalent to considering a time-dependent $\mathbb{R}^{q'}$-valued variable $z(t)$ that satisfies
$$
z^{(n')} = -\sum_{j=0}^{n'-1} A_{j+1} z^{(j)} + \bar B \dot w  ,
$$
where the notation $(.)^{(j)}$ represents the derivative at the order $j$ with respect to time and $Cy = \bar{C} z$.
\end{remark}

\section{Explicit sensitivity formulas for a particle in a harmonic potential}
\label{sec:harmonic}
The examples we consider are similar to Szamel's configuration.
The correlation time is $\tau_c = \sqrt{2} \tau_p$, and $\tau_p$ will take the same range of values given in \cite{Szamel_2017}. 
In this section, we consider the case of a single one-dimensional particle in a harmonic potential.
The real-valued state variable $x$ satisfies 
\begin{equation}
\label{eq:xkx}
\dot x = \xi_0^{-1}(-kx+f) ,
\end{equation}
where $k>0$ and $f$ is a real-valued Gaussian process with PSD \eqref{eq:psd1:part} or \eqref{eq:psd2:part}.
When $f$ corresponds to \eqref{eq:psd1:part}, then $f=Cy$ and $(x,y)$ satisfies (\ref{eq:x}-\ref{eq:y}) with $F(x) = - k x$. Furthermore, the parameters from \eqref{eq:matrices_for_PSD1} are $n=r=1$, $p=q=2$,
$$
A = 
\ell_1 \begin{pmatrix}
1 & -1/2\\
1/2 & 1
\end{pmatrix}, 
\quad B = \ell_1 I_2, 
\quad C = 
\begin{pmatrix}
\sigma_1 & 0
\end{pmatrix},
$$
and the value of $(\sigma_1,\ell_1)$ is given in \eqref{eq:parameters_psd1}. 
Similarly, when $f$ corresponds to \eqref{eq:psd2:part}, $(x,y)$ also satisfies (\ref{eq:x}-\ref{eq:y}) with $F(x) = - k x$, but the parameters from \eqref{eq:matrices_for_PSD2} are $n=1$, $r=2$, $p=1$, $q=2$, 
$$
A = 
\ell \begin{pmatrix}
0 & 1\\
-1 & -2
\end{pmatrix}, 
\quad B = \ell 
\begin{pmatrix}
0\\
1
\end{pmatrix},
\quad C = 
\begin{pmatrix}
\sigma & 0
\end{pmatrix} ,
$$ 
and the value of $(\sigma,\ell)$ is given in \eqref{eq:parameters_psd2}.
Below, we derive explicit formulas for the sensitivities and deduce their time asymptotic.
The goal of this section is to make it clear that the sensitivity indices depend on the form of the spectrum of the driving noise, as we can get closed-form expressions in this simple example.

We perturb the system in two different ways:
\begin{itemize}
\item by a constant force $\lambda_1$, i.e. $- k x$ is replaced by $- k x + \lambda_1$ in \eqref{eq:xkx}. 
We can write $- k x + \lambda_1$ as $F(x) + \lambda_1 \hat F(x)$ where $\hat F (x) = 1$,
\item by a force constant $\lambda_2$, $- kx$ is replaced by $- kx + \lambda_2 x$ in \eqref{eq:xkx}. We can write $- kx + \lambda_2 x$ as $F(x) + \lambda_2 \hat F(x)$ where $\hat F (x) = x$.
\end{itemize}

\begin{figure*}[ht!]
\centering
\includegraphics[scale=1]{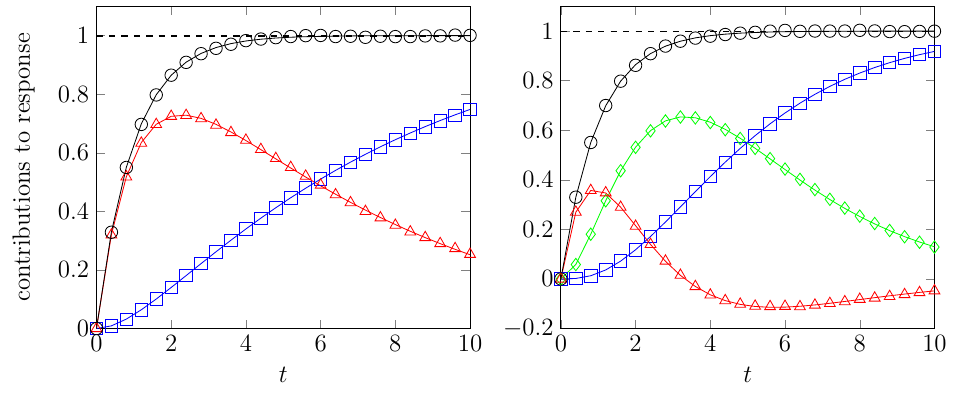}
\caption{First sensitivity index in \eqref{eq:decomp1} and \eqref{eq:decomp2}. 
The solid lines represent the analytical results. The shapes correspond to Monte Carlo simulation results.
Left: PSD \eqref{eq:psd1:part}.
The blue, red, and black lines correspond to $\avg{x(t)p_{0,0}^1(t)}$,$\avg{\dot x(t)p_{1,1}^1(t)}$, and $\frac{\dd}{\dd \lambda_1} \avg{x(t)}_{\lambda_1}$, respectively.
Right: PSD \eqref{eq:psd2:part}. The blue, green, red, and black lines correspond to $\avg{x(t)p_{0,0}^2(t)}$, $\avg{\dot x(t)p_{1,1}^2(t)}$,$\avg{\ddot x(t)p_{2,2}^2(t)}$, and $\frac{\dd}{\dd \lambda_1} \avg{x(t)}_{\lambda_1}$, respectively.} \label{fig:sensitivity_index1}
\end{figure*}

\begin{figure*}[ht!]
\centering
\includegraphics[scale=0.95]{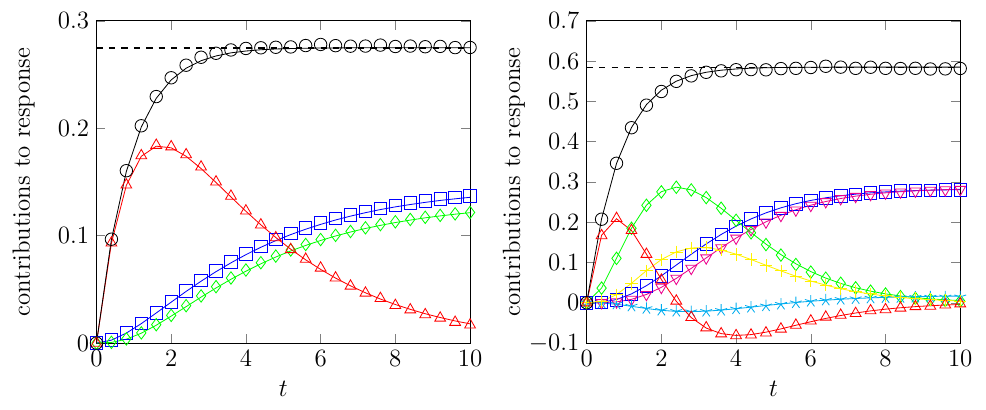}
\caption{
Second sensitivity index in \eqref{eq:decomp3} and \eqref{eq:decomp4}.
The solid lines represent the analytical results. The shapes correspond to Monte Carlo simulation results.
Left: PSD \eqref{eq:psd1:part}. 
The blue, green, red, and black lines correspond to $\avg{x^2(t)p_{0,0}^3(t)}$, 
 $\avg{x^2(t)p_{1,0}^3(t)}$,
 $\avg{\dot{x^2}(t)p_{1,1}^3(t)}$, 
 and 
 $\frac{\dd}{\dd \lambda_2} \avg{x^2(t)}_{\lambda_2}$, respectively.
Right: PSD \eqref{eq:psd2:part}. The blue, magenta, green, cyan, yellow, red, and black lines correspond to $\avg{x^2(t)p_{0,0}^4(t)}$, $\avg{x^2(t)p_{1,0}^4(t)}$,$\avg{\dot{x^2}(t)p_{1,1}^4(t)}$,$\avg{x^2(t)p_{2,0}^4(t)}$,$2\avg{\dot{x^2}(t)p_{2,1}^4(t)}$,$\avg{\ddot{x^2}(t)p_{2,2}^4(t)}$, and $\frac{\dd}{\dd \lambda_2} \avg{x_{\lambda_2}^2(t)}$, respectively. }
\label{fig:sensitivity_index2}
\end{figure*}

Calculations show that the sensitivity of the first (second) moment of $x(t)$ with respect to the second (first) perturbation under the two types of noise \eqref{eq:psd1:part} and \eqref{eq:psd2:part} is zero.  Therefore, below, we focus on the sensitivity of the first (second) moment of $x(t)$ with respect to the first (second) perturbation.  
We have calculated the first sensitivity index $\dd \avg{x(t)}_{\lambda_1}/ \dd \lambda_1$ in two different ways. 
The first one relies on the evaluation of the weighted averages that appear on the right-hand side (rhs) of \eqref{eq:sensitivity1} 
and \eqref{eq:sensitivity2} and summing them up.
The second one, for verification, consists of evaluating the first moment of the perturbed dynamics and then differentiating with respect to the perturbation parameters.
In each case, both results agree
(see Figures \ref{fig:sensitivity_index1}-\ref{fig:sensitivity_index2}).
We have also estimated the weighted averages with Monte Carlo simulations that follow our analytical predictions.
Here, for the two types of noise with PSD \eqref{eq:psd1:part} and \eqref{eq:psd2:part}, the expression is simple: we get 
$\dd \avg{x(t)}_{\lambda_1}/ \dd \lambda_1 |_{\lambda_1=0}
= k^{-1}(1-e^{-k\xi_0^{-1}t}).$
Nonetheless, the two decompositions in terms of weighted averages are different.
For Equation \eqref{eq:sensitivity1} that is related to PSD \eqref{eq:psd1:part}, 
the decomposition is as follows
\begin{equation}
\label{eq:decomp1}
 \frac{\dd}{\dd \lambda_1} \avg{ x(t) }_{\lambda_1} \Big|_{\lambda_1 = 0} = \avg{ x(t) p_{0,0}^1(t) } + \avg{ \dot x(t) p_{1,1}^1(t) } ,
\end{equation}
where the Malliavin weights are $p_{0,0}^1(t) = \sigma_1^{-1} w_1(t) + \sigma_1^{-1} \ell_1^{-1} \omega_1 w_2(t)$, and \\ $p_{1,1}^1(t) = \sigma_1^{-1} \ell_1^{-1} w_1(t)$. 
Note that $p_{1,0}^1(t) = 0$.
Then, for Equation \eqref{eq:sensitivity2} that is related to PSD \eqref{eq:psd2:part} the decomposition is follows
\begin{equation}
\label{eq:decomp2}
    \frac{\dd}{\dd \lambda_1} \avg{ x(t) }_{\lambda_1} \Big|_{\lambda_1 = 0} = \avg{x(t) p_{00}^2(t)} + \avg{\dot x(t) p_{11}^2(t)} + \avg{\ddot x(t) p_{22}^2(t)}   ,
\end{equation}
where the Malliavin weights are
$p_{0,0}^2(t) = \sigma^{-1}\alpha_0 w(t), p_{1,1}^2(t) = \sigma^{-1} \ell^{-1} \alpha_1 w(t)$, and $p_{2,2}^2(t) = \sigma^{-1}\ell^{-2} w(t)$.
 See Appendix \ref{appendix:weights} for details on the weighted averages.
Similarly, we have also calculated the second sensitivity index $\dd \avg{x(t)^2}_{\lambda_2}/ \dd \lambda_2$ in two different ways. 
The counterpart of \eqref{eq:decomp1} is 
\begin{equation}
\label{eq:decomp3}
    \frac{\dd}{\dd \lambda_2} \avg{ x(t)^2 }_{\lambda_2} \Big|_{\lambda_2 = 0} = 
    \avg{x^2(t)p_{0,0}^3(t)} 
    + \avg{x^2(t)p_{1,0}^3(t)} 
    + \avg{\dot{x^2}(t)p_{1,1}^3(t)}
\end{equation}
where 
\begin{align*}
    & p_{0,0}^3(t) = \sigma_1^{-1} \int_0^t x(s) dw_{1}(s) + \sigma_1^{-1} \ell_1^{-1} \omega_1 \int_0^t x(s) dw_{2}(s), \quad 
    p_{1,0}^3(t) = \sigma_1^{-1} \ell_1^{-1} \int_0^t \dot x(s) dw_1(s),
    \\& \quad \text{and} \quad p_{1,1}^3(t) = \sigma_1^{-1} \ell_1^{-1} \int_0^t x(s) dw_1(s).
\end{align*}

The counterpart of \eqref{eq:decomp2} is 
\begin{align}
\nonumber
    \frac{\dd}{\dd \lambda_2} \avg{x^2(t)}_{\lambda_2}\Big|_{\lambda_2 = 0}  = & \avg{x^2(t)p_{0,0}^4(t)} + \avg{x^2(t)p_{1,0}^4(t)} + \avg{x^2(t)p_{2,0}^4(t)}\\
    & + \avg{\dot {x^2}(t)p_{1,1}^4(t)}
    + 2\avg{\dot {x^2}(t)p_{2,1}^4(t)}
    + \avg{\ddot {x^2}(t)p_{2,2}^4(t)}  ,
    \label{eq:decomp4}
\end{align}
where 
\begin{align*}
& p_{0,0}^4(t) = \sigma^{-1} \alpha_0 \int_0^t x(s) dw(s),
\: \: 
p_{1,1}^4(t) = \sigma^{-1} \alpha_1 \ell^{-1} \int_0^t x(s) dw(s), 
\: \: \\
& p_{2,2}^4(t) = \sigma^{-1} \ell^{-2} \int_0^t x(s) dw(s), 
\: \:
p_{1,0}^4(t) = \sigma^{-1} \alpha_1
\ell^{-1} \int_0^t \dot x(s) dw(s),
\: \: \\
& p_{2,1}^4(t) = \sigma^{-1} \ell^{-2} \int_0^t \dot x(s) dw(s), 
\: \:  
\mathrm{and} \: \: p_{2,0}^4(t) = \sigma^{-1}\ell^{-2} \int_0^t \ddot x(s) dw(s).
\end{align*}
In each case, both analytical results agree.
We represent our analytical finding in Figure \ref{fig:sensitivity_index2} including Monte Carlo simulations.
Among the two sensitivity indices above, only the second depends on the correlation time.
This is consistent with \cite{Szamel_2017} for the case of a simple OU process.
We also obtain the long-time asymptotic.
In all cases, for the first sensitivity index, we have 
$$ 
\lim \limits_{t \to \infty} \frac{\dd \avg{x(t)}_{\lambda_1}}{\dd \lambda_1}\Big|_{\lambda_1=0}
= \frac{1}{k}.
$$
In contrast, the long-time asymptotic of 
the second sensitivity index depends on the type of noise. For \eqref{eq:psd1:part},
$$ 
\lim \limits_{t \to \infty} \frac{\dd \avg{x(t)^2}_{\lambda_2}}{ \dd \lambda_2} \Big|_{\lambda_2 = 0} 
= \frac{\sigma_1^2}{2k^2\xi_0} \frac{(2k(\xi_0\ell_1)^{-1}+1)(k(\xi_0\ell_1)^{-1}+1)^2+(\frac{\omega}{\ell_1})^2}{((k(\xi_0\ell_1)^{-1}+1)^2+(\frac{\omega}{\ell_1})^2)^2}.
$$
For \eqref{eq:psd2:part},
$$ 
\lim \limits_{t \to \infty} \frac{\dd \avg{x(t)^2}_{\lambda_2}}{ \dd \lambda_2} \Big|_{\lambda_2 = 0} 
= \frac{\sigma^2}{4k^2\xi_0} \frac{k/(\xi_0\ell)+2}{(k/(\xi_0\ell)+1)^2} + \frac{\sigma^2}{4k\xi_0^2\ell} \frac{ (k/(\xi_0\ell)+2)^2-1}{(k/(\xi_0\ell)+1)^4}.
$$
Such a difference in behavior
can also be observed in the steady state/long-time asymptotic variance of the response. 
For \eqref{eq:psd1:part},
$$ 
\lim \limits_{t \to \infty} \avg{x(t)^2}
=  \frac{\ell_1(\frac{k}{\xi_0}+\ell_1)\sigma_1^2}{2k\xi_0((\frac{k}{\xi_0}+\ell_1)^2+\omega_1^2)},
$$
For \eqref{eq:psd2:part},
$$ 
\lim \limits_{t \to \infty}  \avg{x(t)^2}
= \frac{\sigma^2}{4k\xi_0} \frac{2+k \xi_0^{-1}\ell^{-1}}{(1+k \xi_0^{-1}\ell^{-1})^2}.
$$
When $\tau_c \to 0$, in the first two equations above both rhs converge towards 
$T_{\textup{eff}}k^{-2}$ 
and in
the last two equations above both rhs converge towards $T_{\textup{eff}}k^{-1}$.
That is consistent with the long-time asymptotic response of a thermal Brownian particle in a harmonic potential.


\section{Sensitivity analysis for an interacting particle system in a screened Coulomb potential}
\label{sec:ips}
We consider $N$ particles in three dimensions interacting via a repulsive
screened Coulomb potential, $\forall r > 0, \: V (r) = A_V \exp (-\kappa (r - \sigma_V)) /r$. 
We use the same parameters as Szamel in \cite{Szamel_2017}:
$N = 1372$, $A_V = 475 T_{\textup{eff}}^{\textup{sp}} \sigma_V$ and $\kappa \sigma_V = 24$, at number density $N \sigma_V^2/V = 0.51$.
The state variable 
$x = \{x_{i\alpha}\}_{1 \leq i \leq N, \: 1 \leq \alpha \leq 3 }$ (stacked in one vector of size $n=3N$) is $n=4116$ dimensional and satisfies the equation
$$
\dot x = \xi_0^{-1}(F(x)+f),
$$
where the unperturbed force for each particle $1 \leq i \leq N$ can be expressed as 
$$
F_{i \alpha}(x) = - \partial_{x_{i \alpha}} \left ( \sum \limits_{j \neq i} V(\| x_i - x_j \|)\right )\\
 = - \sum \limits_{j \neq i} (x_{i \alpha} - x_{j \alpha})\frac{V'(\| x_i - x_j \|)}{\| x_i - x_j \|}, \: \: \alpha = 1,2,3,
$$
and the colored noise $f$ is a $n$-dimensional zero-mean SGP with PSD \eqref{eq:psd1:part} or \eqref{eq:psd2:part}. We use the notation $\| u \| = \sqrt{u_1^1 + u_2^2 + u_3^2}$.
When $f$ corresponds to \eqref{eq:psd1:part}, then $f=Cy$ and $(x,y)$ satisfies Equation (\ref{eq:x}-\ref{eq:y}) with matrices $A, B,C$ given in \eqref{eq:matrices_for_PSD1} with $n=4116$, $r=1$, $p=q= 2rn =8232$,
and the value of $(\sigma_1,\ell_1)$ is given in \eqref{eq:parameters_psd1}.
Similarly, when $f$ corresponds to \eqref{eq:psd2:part}, $(x,y)$ satisfies also Equation (\ref{eq:x}-\ref{eq:y}) with matrices $A, B,C$ given in \eqref{eq:matrices_for_PSD2} with $n=4116$, $r=2$, $p= n = 4116$, $q = 8232$,
and the value of $(\sigma,\ell)$ is given in \eqref{eq:parameters_psd2}.
Next, we evaluate the time-dependent mobility of a single particle in the system of interacting particles driven by the $n$-dimensional noise $f$ (self-propulsion) following \eqref{eq:psd1:part} and \eqref{eq:psd2:part}.  
Since explicit formulas are unavailable, we perform computations using our sensitivity formula \eqref{eq:sensitivity1}-\eqref{eq:weight1}, \eqref{eq:sensitivity2}-\eqref{eq:secondgeneralization}, and unperturbed simulations. 
We fix a particle $1 \leq i \leq N$ and a direction $1 \leq \alpha \leq 3$ and we perturb the system at $t=0^+$ as follows:  
$F_i(x)$ is replaced by $F_i(x) + \lambda \boldsymbol{e}_\alpha$ 
where $\boldsymbol{e}_\alpha$ is the unit vector in the direction $\alpha$.
Under the influence of the additional constant force, the component $\alpha$ of the particle $i$ moves according to the mobility function $\chi(t)$ as follows $\langle x_{i\alpha}(t) \rangle_\lambda = \chi(t) \lambda + o(\lambda)$ and thus $\dd \avg{x(t)}_{\lambda}/\dd \lambda_{| \lambda = 0} = \chi(t)$. Moreover,  $\lim \limits_{t \to \infty} \chi(t)/t = \mu$, the mobility coefficient.
In principle, to estimate $\chi(t)$ we can compute the sensitivity index above with unperturbed averages 
$\avg{x_{i\alpha}(t) p_{00;i\alpha}^1(t)}
+ \avg{\dot x_{i\alpha}(t) p_{11;i\alpha}^1(t)}$ or $\avg{x_{i\alpha}(t) p_{00;i\alpha}^2(t)} + \avg{\dot x_{i\alpha}(t) p_{11;i\alpha}^2(t)} + \avg{\ddot x_{i\alpha}(t) p_{22;i\alpha}^2(t)}$
when the driving noise corresponds to \eqref{eq:psd1:part} or \eqref{eq:psd2:part}, respectively.  Here the corresponding weight functions evolve according to the following equations of motion:
\begin{itemize}
\item 
For \eqref{eq:psd1:part}, 
$$
p_{00;i\alpha}^1(t) = 
\sigma_1^{-1} (w_{i\alpha})_1(t) + \sigma_1^{-1} \ell_1^{-1} \omega_1 (w_{i\alpha})_2(t),
\: \mbox{ and } \:
p_{11;i\alpha}^1(t) = \sigma_1^{-1} \ell_1^{-1} (w_{i\alpha})_1(t).
$$
\item
For \eqref{eq:psd2:part}, 
$$
p_{00;i\alpha}^2(t) =\sigma^{-1}\alpha_0 w_{i\alpha}(t), \:
p_{11;i\alpha}^2(t) =\sigma^{-1} \ell^{-1} \alpha_1 w_{i\alpha}(t), \: \mbox{ and } \:
p_{22;i\alpha}^2(t) =
\sigma^{-1}\ell^{-2} w_{i\alpha}(t).
$$
\end{itemize}

\begin{figure*}[b!]
\centering
\includegraphics[scale=0.5]{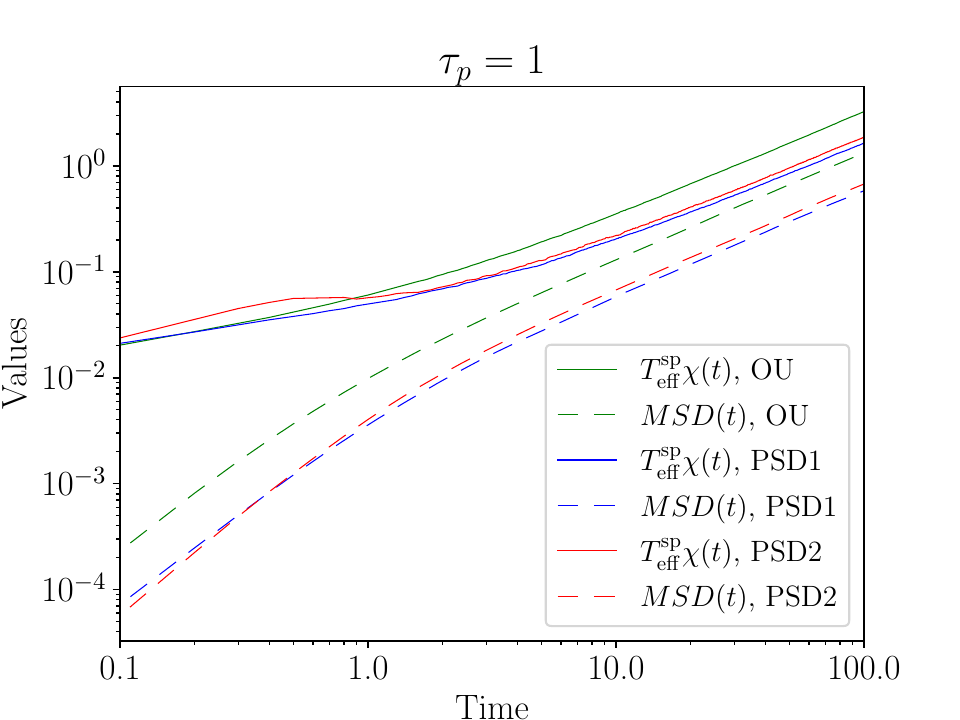}
\includegraphics[scale=0.5]{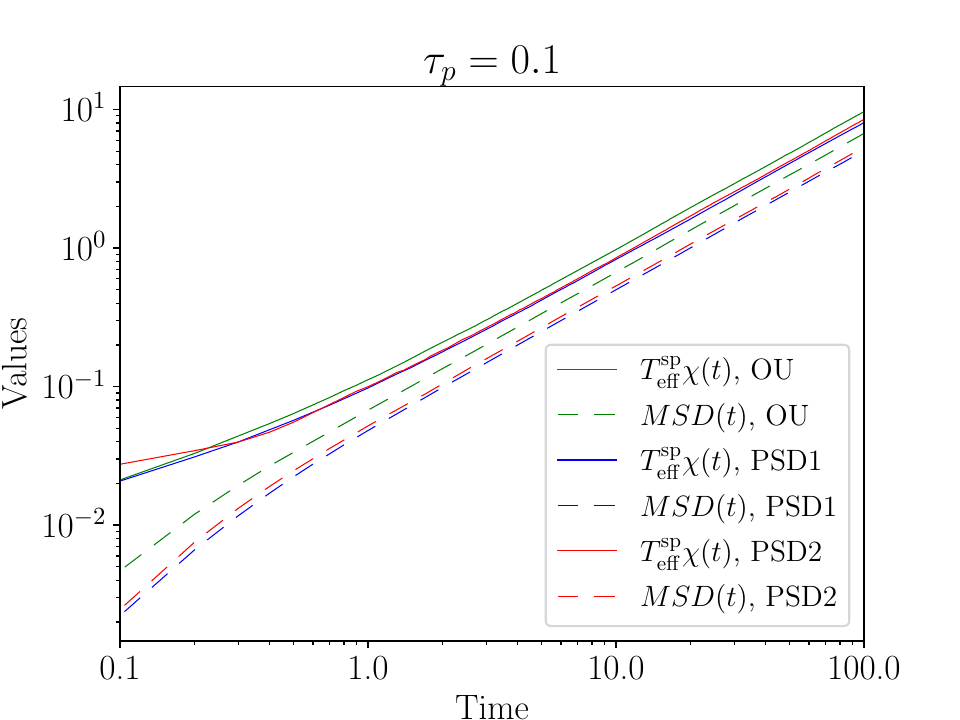}
\includegraphics[scale=0.5]{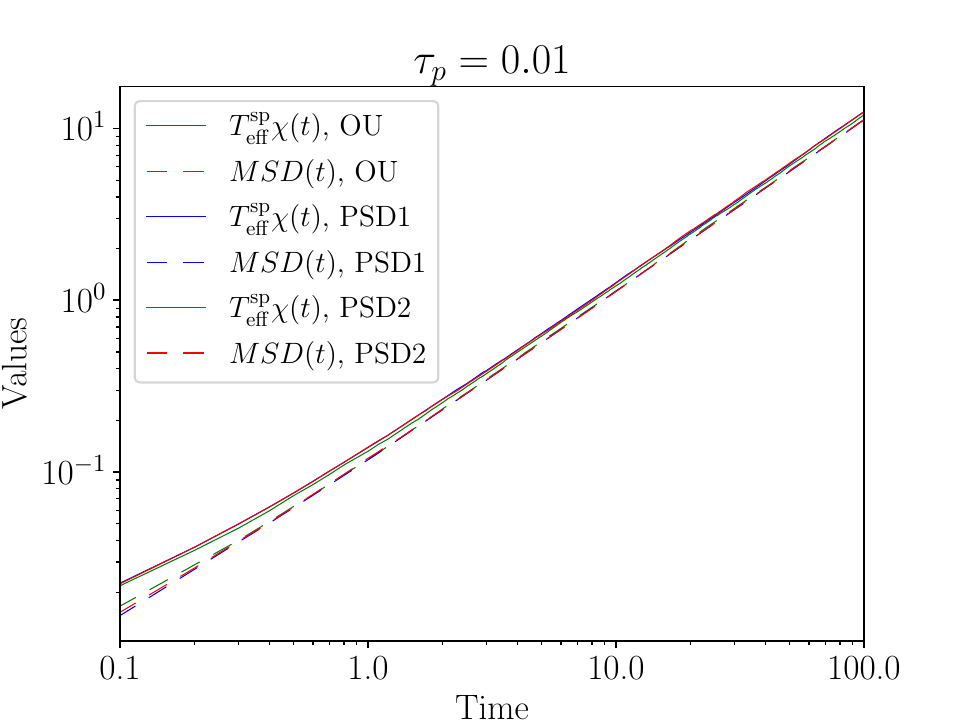}
\includegraphics[scale=0.5]{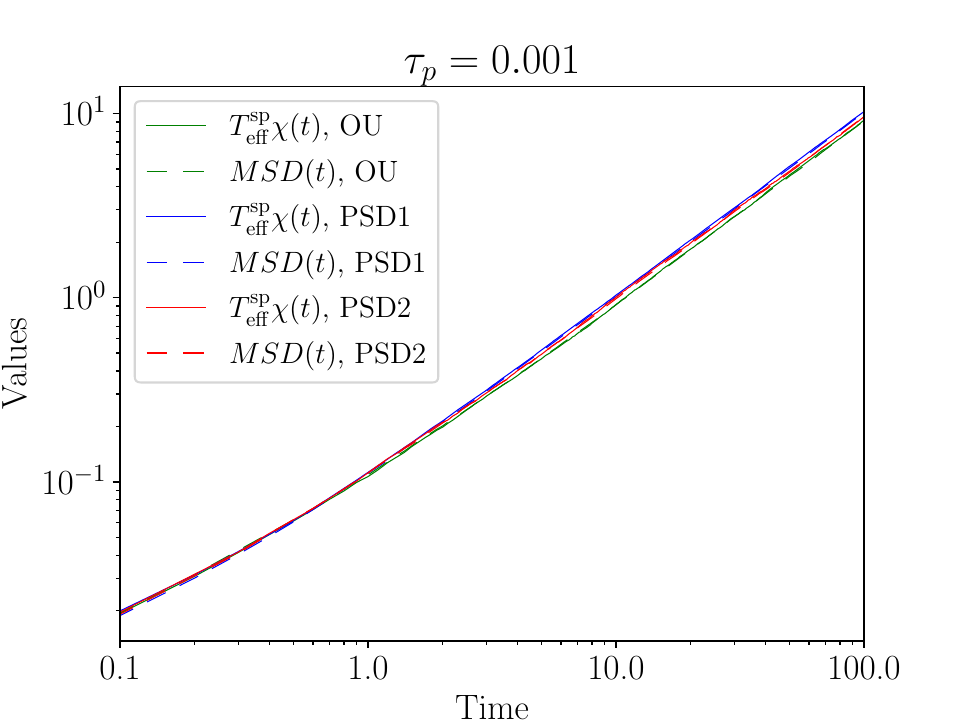}
\caption{Mobility function, $T_\text{eff}^\text{sp} \chi(t)$, (solid) and mean-square displacement, $MSD(t)$, (dashed) as a function of time, for $\tau_p = 1 \text{(top left)}, 0.1 \text{(top right)}, 0.01 \text{(bottom left)}, 0.001 \text{(bottom right)}$.
Superposition of the results for OU (green), PSD1 (blue), PSD2 (red) with spectra \eqref{eq:psd0}, \eqref{eq:psd1:part}, and  \eqref{eq:psd2:part}, respectively. 
}
\label{figMSDCHI}
\end{figure*}

\begin{figure}[ht!]
\centering
\includegraphics[scale=1]{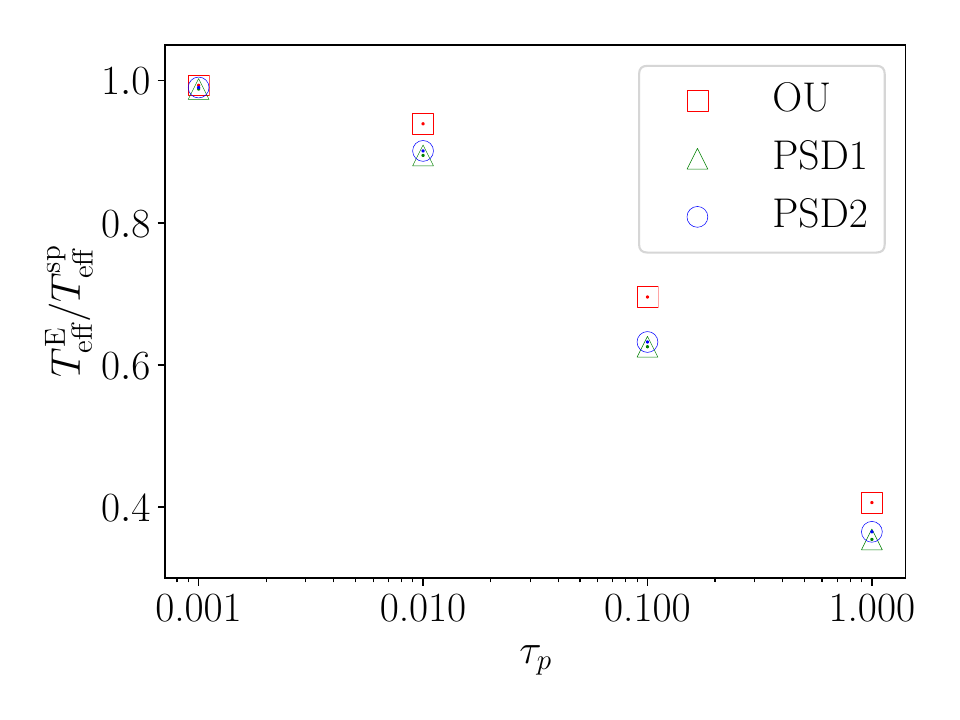}
\caption{$T_{\rm eff}^E/T_{\rm eff}^{ \rm sp}$, as a function of
persistence time, $\tau_p$.
Superposition of the results for OU (green), PSD1 (blue), PSD2 (red) with spectra \eqref{eq:psd0}, \eqref{eq:psd1:part}, and  \eqref{eq:psd2:part}, respectively. 
}
\label{fig:EffT}
\end{figure}

For both cases, since the particles are indistinguishable and exhibit isotropic behavior, the mobility function $\chi(t)$ does not depend on $(i,\alpha)$. Then, as inspired by Szamel's computing strategy, we can use all the particles in all Cartesian directions and average over the origins of time.
The mobility function $\chi(t)$ can thus be estimated by 
\begin{equation}
\label{eq:numerator}
\frac{1}{3 N N_o}
\sum \limits_{i=1}^N \sum \limits_{\alpha=1}^3
\sum \limits_{\mathfrak{o}=1}^{N_o} \zeta(i,\alpha,\mathfrak{o},t),
\end{equation}
where in the case of \eqref{eq:psd1:part},
\begin{align*}
\zeta(i,\alpha,\mathfrak{o},t) =  
& \avg{x_{i\alpha}(t+t_{\mathfrak{o}}) (p_{00;i\alpha}^1(t+t_{\mathfrak{o}}) - p_{00;i\alpha}^1(t_{\mathfrak{o}}))} \\
& + \avg{\dot x_{i\alpha}(t+t_{\mathfrak{o}}) (p_{11;i\alpha}^1(t+t_{\mathfrak{o}}) - p_{11;i\alpha}^1(t_{\mathfrak{o}})})
\end{align*}
and in the case of \eqref{eq:psd2:part},
\begin{align*}
\zeta(i,\alpha,\mathfrak{o},t) = & \avg{x_{i\alpha}(t+t_{\mathfrak{o}}) (p_{00;i\alpha}^2(t+t_{\mathfrak{o}})-p_{00;i\alpha}^2(t_{\mathfrak{o}})} \\
& + \avg{\dot x_{i\alpha}(t+t_{\mathfrak{o}}) (p_{11;i\alpha}^2(t+t_{\mathfrak{o}})-p_{11;i\alpha}^2(t+t_{\mathfrak{o}})}\\ 
&  + \avg{\ddot x_{i\alpha}(t+t_{\mathfrak{o}}) (p_{22;i\alpha}^2(t+t_{\mathfrak{o}})-p_{22;i\alpha}^2(t_{\mathfrak{o}})}.
\end{align*}
Finally,
in both cases 
we can compute an estimation of the mean-square displacement (MSD) defined by
\begin{equation}
\label{eq:denominator}
{\rm MSD}(t)=
\avg{\| x_1(t) - x_1(0)\|^2} = 
\frac{1}{N} \sum \limits_{i=1}^N
\avg{\| x_i(t) - x_i(0)\|^2},
\end{equation}
which grows as $6D_{\rm sd}t$, with $D_{\rm sd}$ being the self-diffusion coefficient.
The Einstein relation can be used to define the Einstein effective temperature $T_{\rm eff}^E = D_{\rm sp} / \mu$ \cite{Szamel_2017}.
It can then be estimated by the ratio of ${\rm MSD}(t)/6$ over $\chi(t)$ for $t$ large. For four distinct persistence times, we plot, in Figure \ref{figMSDCHI}, the mobility function and mean-square displacement for the systems corresponding to the three noises that we consider.
It is known that the Einstein effective temperature can be different from $T_{\rm eff}^{\rm sp}$, except in the small correlation time limit.
In Figure \ref{fig:EffT}, we show that the ratio $T_{\rm eff}^E /T_{\rm eff}^{\rm sp}$ decreases when the correlation time increases and that the decay also depends on the form of the PSD. The value is, therefore, a complicated function of the PSD and not only a function of the amplitude and correlation time of the colored noise.

\section{Conclusion}
In this paper, we have proposed a method of sensitivity analysis for particle systems under colored noise. We have obtained original formulas, generalizing those obtained by Szamel, for the sensitivity indices of the IPS statistics mentioned above with respect to the perturbation of the drift coefficient.
The type of colored noise we considered falls into a class of stationary Gaussian processes with zero mean and whose spectra can be general. Both our analytical and numerical calculations show that the sensitivity indices do not depend only on the effective parameters: the amplitude and coherence time of the noise. Instead, it is also dependent on the structure of the spectrum. Moreover, the method that we have developed can be applied beyond the quantities studied in this paper. For example, it can be used to study the structure of correlations in the particle system, as done in \cite{maggi2022critical} (case of a simple OU). Furthermore, our method can prove effective for gradient calculation in optimization procedures, as was done in \cite{das2023nonequilibrium} (case of white noise).

\section{Appendix: Proofs and details of weight formulas}
\subsection{Derivation of our first main result}
\label{appendix:proof_formula_1}
We extend the approach in \cite{Szamel_2017} to compute the derivative of $\langle \Phi(x(t)) \rangle_\lambda$  with respect to $\lambda$ at $\lambda = 0$.

\noindent \textbf{Step 1}\\ 
We start by discretizing the equations of motion (\ref{eq:x}-\ref{eq:y}), for fixed $N_t$, over time intervals of length $\Delta t = t/N_t$ to obtain the following limit of integrals
\begin{align*}
    \langle \Phi(x(t)) \rangle_\lambda = \lim_{N_t \rightarrow \infty} \int \Phi(x_{N_t}) \prod_{i=2}^{N_t} P_{\lambda,i} P_1 P^{ss},
\end{align*}
where $P_{\lambda,i}=P_\lambda(x_i,y_i|x_{i-1},y_{i-1})$, for $i=2,...,{N_t}$, are the transition densities of the perturbed system, $P_1=P(x_1,y_1|x_{0},y_{0})$ is the transition density of the unperturbed system, and $P^{ss}=P^{ss}(x_0,y_0)$ is the steady state distribution for the unperturbed system.
Here, the transition densities, for $i=2,...,N_t$, are given by
\begin{align*}
    P_\lambda(x_i,y_i|x_{i-1},y_{i-1}) = \delta_{\mathcal{X}_{\lambda,i}} g(y_i,(I_q - \Delta t A) y_{i-1},\Delta t B B^T),
\end{align*}
where $\mathcal{X}_{\lambda,i} = x_i - x_{i-1} - \Delta t \xi_0^{-1} \left ( F_\lambda(x_{i-1}) + C y_{i-1} \right )$ and $g(y,\mu,\Sigma)$ is the density of the multidimensional Gaussian $\mathcal{N}(\mu,\Sigma)$. We use the notation $I_q$ for the identity matrix of size $q \times q$.\\ 

\noindent \textbf{Step 2}\\
To compute the derivative $\frac{\dd}{ \dd \lambda} \langle \Phi(x(t)) \rangle_\lambda $, we interchange the limit and integral on the rhs of the equation, followed by applying the product rule to obtain the following
\begin{align*}
    \frac{\dd}{\dd \lambda} \langle \Phi(x(t)) \rangle_\lambda \Big|_{\lambda=0}= \lim_{{N_t} \rightarrow \infty} \sum_{i=2}^{N_t} \int \Phi(x_{N_t}) P_{\lambda=0,i}' \left ( \prod_{j=2, j \neq i}^{N_t} P_j \right ) P_1 P^{ss}.
\end{align*}
Using the assumption that $\hat F=C E$ and the chain rule, it can be shown that
\begin{align*}
    P_{\lambda=0,i}' = E (x_{i-1}) \cdot \nabla_{y_{i-1}}  \left ( \delta_{\mathcal{X}_{\lambda,i}} \right ) g(y_i,(I_q - \Delta t A) y_{i-1}, \Delta t BB^T).
\end{align*}

\noindent \textbf{Step 3}\\
By applying integration by parts and rearranging, we have the following
\begin{align*}
    & \frac{\dd}{\dd \lambda} \langle \Phi(x(t)) \rangle_\lambda \Big|_{\lambda=0}=\\ 
    & \lim_{{N_t} \rightarrow \infty} \int \Phi(x_{N_t}) \prod_{j=1}^{N_t} P_j P^{ss}
    \Bigg\{ - \bigg[ (I_q - \Delta t A) E (x_{{N_t}-1}) \cdot (\Delta t BB^T)^{-1}
    \\& \qquad\qquad\qquad\qquad\qquad\qquad\qquad \times (y_{N_t} - (I_q - \Delta t A)y_{N_t-1}) \bigg]
    \\& + \sum_{i=1}^{{N_t}-1} \bigg[ \frac{(E (x_i) - E (x_{i-1}))}{\Delta t} \cdot (\Delta t BB^T)^{-1} (y_i - (I_q - \Delta t A)y_{i-1}) \Delta t \bigg]
    \\& + \sum_{i=2}^{{N_t}-1} \bigg[ \Delta t  A E (x_{i-1}) \cdot (\Delta t BB^T)^{-1} (y_i - (I_q - \Delta t A)y_{i-1}) \bigg]
    \\& + \bigg[ E (x_0) \cdot (\Delta t BB^T)^{-1} (y_1 - (I_q - \Delta t A)y_0) \bigg] \Bigg\},
\end{align*}
Note that this equation generalizes Equation A.3 in \cite{Szamel_2017}.\\

\noindent \textbf{Step 4}\\
It can be shown that an equivalent identity to Equation A.4 in \cite{Szamel_2017} can be derived from the steady state property, which results in the following
\begin{align*}
    & \int \Bigg[ \Phi(x_{N_t}) E(x_{0}) \cdot (\Delta t BB^T)^{-1} (y_1 - (I_q - \Delta t A)y_0) \cdots\\
    & \cdots - \Phi(x_{N_t}) E(x_{{N_t}-1}) \cdot (\Delta t BB^T)^{-1} (y_{N_t} - (I_q - \Delta t A)y_{N_t-1}) \Bigg] \prod_{j=1}^{N_t} P_j P^{ss}\\
    & = \int \frac{\Phi(x_{N_t})-\Phi(x_{{N_t}-1})}{\Delta t}
    \Bigg[ \sum_{i=1}^{{N_t}-1} E(x_{i-1}) \cdot (BB^T)^{-1} (y_i - (I_q - \Delta t A)y_{i-1}) \Bigg] \prod_{j=1}^{N_t} P_j P^{ss}.
\end{align*}

\noindent \textbf{Step 5}\\
The equation can now be written in the form of finite differences as below
\begin{align*}
    & \frac{\dd}{\dd \lambda}  \langle \Phi(x(t)) \rangle_\lambda\Big|_{\lambda=0}
    \\
    & = \lim_{{N_t} \rightarrow \infty} \int \Bigg\{ \Phi(x_{N_t}) \Bigg[ \sum_{i=2}^{N_t} A E(x_{i-1}) \cdot (BB^T)^{-1} (y_i - (I_q - \Delta t A)y_{i-1}) \Bigg]
    \\
    & \quad + \Phi(x_{N_t}) \Bigg[ \sum_{i=1}^{{N_t}-1} \nabla  E(x_{i-1}) \frac{x_i-x_{i-1}}{\Delta t} \cdot (BB^T)^{-1}(y_{i}-(I_q - \Delta t A)y_{i-1}) \Bigg]
    \\
    & \quad  + \frac{\Phi(x_{N_t})-\Phi(x_{{N_t}-1})}{\Delta t}
    \Bigg[ \sum_{i=1}^{{N_t}-1} E(x_{i-1}) \cdot (BB^T)^{-1} (y_i - (I_q - \Delta t A)y_{i-1}) \Bigg] \Bigg\} \prod_{j=1}^{N_t} P_j P^{ss}.
\end{align*}
This corresponds to Equation A.5 in \cite{Szamel_2017}.\\

\noindent \textbf{Step 6}\\
Assuming that limits can be taken, we obtain 
\eqref{eq:sensitivity1} and 
\eqref{eq:weight1}.
\subsection{Derivation of our second main result}
\label{app:proof2}
Our proof is composed of several steps.
Starting with the third step, we provide the proof in details for the case $n'=2$. We give the main formula for the general case $n'\geq 1$.\\

\noindent \textbf{Step 1}\\
We start the perturbation at $x_{n'+1}$.
Here, the transition densities, for $i=n'+1,...,N_t$, are given by
\begin{align*}
    & P_\lambda(x_i,y_i|x_{i-1},y_{i-1}) \\
    & \qquad = \delta_{\mathcal{X}_{\lambda,i}} \left ( \prod \limits_{k=1}^{n'-1} \delta_{\mathcal{Y}_{i,k}} \right ) \: g \left (y_{i,n'}, 
(I_{q'} - \Delta t A_{n'}) y_{i-1,n'} - \Delta t \sum \limits_{k=1}^{n'-1} A_k y_{i-1,k}
, \Delta t \bar B \bar B^T 
\right ),
\end{align*}
where $\mathcal{X}_{\lambda,i}$ is defined as above and $\mathcal{Y}_{i,k} = \{ y_{i,k} - y_{i-1,k} - \Delta t y_{i-1,k+1} \}$ and $g(y,\mu,\Sigma)$ is the density of the multidimensional Gaussian $\mathcal{N}(\mu,\Sigma)$. We use the notation $I_{q'}$ for the identity matrix of size $q' \times q'$.
The unperturbed transition densities for $i=1,...,n'$ are given by the same formula except that $\lambda = 0$. Below, we use the short notation $P_{\lambda, i}$, and $P_i$ for the perturbed and unperturbed transition densities respectively.\\

\noindent \textbf{Step 2}\\
To compute the derivative $\frac{d}{d
\lambda} \langle \Phi(x(t)) \rangle_\lambda$ at $\lambda = 0$, we interchange the limit and integral on the rhs of the equation, followed by applying the product rule to obtain the following
\begin{align*}
    \frac{d}{d \lambda} \langle \Phi(x(t)) \rangle_\lambda \Big|_{\lambda=0} = \lim_{N_t \rightarrow \infty} \sum_{i=n'+1}^{N_t} \int \Phi(x_{N_t}) P_{\lambda=0,i}' \left ( \prod_{j=1, j \neq i}^{N_t} P_j \right ) P^{ss}.
\end{align*}
We will use the short notation $\mathcal{X}_i$ for $\mathcal{X}_{0,i}$. Using the assumption that $\hat F=\bar C \bar E$ and the chain rule, it can be shown that
\begin{align*}
    & P_{\lambda=0,i}' = \\
    & \bar E (x_{i-1}) \cdot \nabla_{y_{i-1,1}}  \left ( \delta_{\mathcal{X}_i} \right ) \left ( \prod \limits_{k=1}^{n'-1} \delta_{\mathcal{Y}_{i,k}} \right ) \: \underbrace{g \left (y_{i,n'}, (I_{q'}-\Delta t A_{n'}) y_{i-1,n'} 
    - \Delta t \sum_{k=1}^{n'-1} A_k y_{i-1,k}, \Delta t \bar B \bar B^T \right )}_{g_i :=}.
\end{align*}
Here we recall that both $\bar E (x_{i-1})$ and  $\nabla_{y_{i-1,1}}(\delta_{\mathcal{X}_{\lambda,i}})$ take values in $\mathbb{R}^{q'}$. Below we use the notation 
\begin{equation}
\label{eq:Ti}
\mathcal{T}_i \triangleq \int \Phi(x_{N_t}) P_{\lambda=0,i}' \left ( \prod_{j=1, j \neq i}^{N_t} P_j \right ) P^{ss}.
\end{equation}
\subsubsection{Case $n'=2$ in detail}
In this subsection we consider the case $n'=2$.
In order to alleviate the notation, we denote $\nabla_u(.)$ as 
$(.)_u$. We reformulate $\mathcal{T}_i$ using several integration by parts (IBPs). We use the notation 
$$
\zeta 
= \Phi(x_{N_t})
\left ( \prod_{j=1}^{N_t} P_j \right ) P^{ss} \: \mbox{ and } \:
\bar E_i = \bar E(x_i),
$$ 
and then
$\mathcal{T}_i$
can be expressed as 
\begin{align}
\label{eq:Ti_reform1}
\mathcal{T}_i = & \frac{1}{\Delta t}  \int \zeta \bar E_{i-1} 
\left ( \left ( \log (g_{i-1} g_{i}) \right)_{y_{i-1,2}} 
- \left ( \log (g_{i-2} g_{i-1}) \right )_{y_{i-2,2}} \right ) \nonumber \\
& + \int \zeta \bar E_{i-1} \: \left ( -\log (g_i) \right )_{y_{i-1,1}}.
\end{align}
The transition from  \eqref{eq:Ti} to \eqref{eq:Ti_reform1} is given in Lemma \ref{appendix:lemma_3132}. Then calculations yield
\begin{align*}
& \left ( -\log (g_i) \right )_{y_{i-1,1}} = 
A_1^T (\bar B \bar B^T)^{-1} \theta_i,\\
& \left ( \log (g_{i-1} g_i) \right )_{y_{i-1,2}}
=
\frac{1}{\Delta t} \left ( \bar B \bar B^T \right )^{-1} \left ( \theta_i-\theta_{i-1}\right ) - A_2^T (\bar B \bar B^T)^{-1} \theta_i,\\
& \left ( \log (g_{i-1} g_{i}) \right )_{y_{i-1,2}} - \left ( \log (g_{i-2} g_{i-1}) \right )_{y_{i-2,2}} \\
& \qquad =
\frac{1}{\Delta t} \left ( \bar B \bar B^T \right )^{-1} \left ( \theta_i-2\theta_{i-1} + \theta_{i-2} \right ) - A_2^T (\bar B \bar B^T)^{-1} (\theta_i-\theta_{i-1})
\end{align*}
where 
$$
\theta_i \triangleq 
y_{i,2} - (I_{q'}-\Delta t A_2) y_{i-1,2}
+ \Delta t A_1y_{i-1,1}.  
$$
Therefore 
\begin{align}
\label{eq:Ti_reform2}
\mathcal{T}_i =
& \int \zeta \bar E_{i-1}
(\bar B \bar B^T)^{-1}
\left ( \frac{\theta_i - 2 \theta_{i-1} + \theta_{i-2}}{\Delta t^2}
\right )
- 
\int \zeta \bar E_{i-1}
A_2^T (\bar B \bar B^T)^{-1}
\left ( \frac{\theta_i - \theta_{i-1}}{\Delta t}
\right )\\
& + \int \zeta \bar E_{i-1}
A_1^T (\bar B \bar B^T)^{-1}
\theta_i. \nonumber
\end{align}
Collecting all the terms, we obtain 
$$
\sum_{i=3}^{N_t} \mathcal{T}_i =
S_3
+S_2
+ S_1
$$
where 
$$
S_1 = \sum_{i=3}^{N_t} \int \zeta \bar E_{i-1}
A_1^T (\bar B \bar B^T)^{-1}
\theta_i,
$$
\begin{align*}
S_2 = & 
- \sum_{i=3}^{N_t} \int \zeta \bar E_{i-1}
A_2^T (\bar B \bar B^T)^{-1}
\left ( \frac{\theta_i - \theta_{i-1}}{\Delta t} \right )\\
= &
 \sum_{i=2}^{N_t-1} \int \zeta \left ( \frac{\bar E_i - \bar E_{i-1}}{\Delta t} \right )
A_2^T (\bar B \bar B^T)^{-1}
\theta_i\\  
& + \frac{1}{\Delta t} \int \zeta \bar E_1 A_2^T (\bar B \bar B^T)^{-1} \theta_2
- \frac{1}{\Delta t} \int \zeta \bar E_{N_t-1}
A_2^T (\bar B \bar B^T)^{-1} \theta_{N_t},
\end{align*}
and 
\begin{align*}
S_3 = &
\sum_{i=3}^{N_t} \int \zeta \bar E_{i-1}
(\bar B \bar B^T)^{-1}
\left ( \frac{\theta_i - 2 \theta_{i-1} + \theta_{i-2}}{\Delta t^2}
\right )\\
= &
\sum_{i=3}^{N_t-2} 
\int \zeta \left ( \frac{\bar E_{i+1}- 2 \bar E_i + \bar E_{i-1}}{\Delta t^2}  \right )
(\bar B \bar B^T)^{-1} \theta_i
\\
& + \int \zeta \bar E_{N_t-2} (\bar B \bar B^T)^{-1}  \left ( \frac{\theta_{N_t-1}}{\Delta t^2} \right )
+ \int \zeta \bar E_{N_t-1} (\bar B \bar B^T)^{-1}  \left ( \frac{\theta_{N_t}}{\Delta t^2}
\right ) \\
& - 2 \left ( \int \zeta \bar E_2 (\bar B \bar B^T)^{-1}  \left ( \frac{\theta_2}{\Delta t^2} \right )
+ \int \zeta \bar E_{N_t-1} (\bar B \bar B^T)^{-1}  \left ( \frac{\theta_{N_t-1}}{\Delta t^2}
\right ) \right )\\
& + \int \zeta \bar E_2 (\bar B \bar B^T)^{-1}  \left ( \frac{\theta_1}{\Delta t^2} \right )
+ \int \zeta \bar E_3 (\bar B \bar B^T)^{-1}  \left ( \frac{\theta_2}{\Delta t^2}
\right ).
\end{align*}
Therefore, we have obtained

\begin{align}
\label{eq:sumTi}
\sum_{i=3}^{N_t} \mathcal{T}_i
= & \sum_{i=3}^{N_t} \int \zeta \bar E_{i-1}
A_1^T (\bar B \bar B^T)^{-1}
\theta_i
+ \sum_{i=2}^{N_t-1} \int \zeta \left ( \frac{\bar E_i - \bar E_{i-1}}{\Delta t} \right ) A_2^T (\bar B \bar B^T)^{-1}
\theta_i \\
& + \sum_{i=1}^{N_t-2} \int \zeta \left ( \frac{\bar E_{i+1}- 2 \bar E_i + \bar E_{i-1}}{\Delta t^2}  \right )
(\bar B \bar B^T)^{-1} \theta_i +   \mathcal{B}_1 +\mathcal{B}_2 , \nonumber
\end{align}
with
\begin{align}
    \mathcal{B}_1 =& \frac{1}{\Delta t} \int \zeta \bar E_1 A_2^T (\bar B \bar B^T)^{-1} \theta_2
- \frac{1}{\Delta t} \int \zeta \bar E_{N_t-1}
A_2^T (\bar B \bar B^T)^{-1} \theta_{N_t} , \nonumber\\
    \mathcal{B}_2 =&  \int \zeta (\bar E_{N_t-2} - 2 \bar E_{N_t-1}) (\bar B \bar B^T)^{-1}  \left ( \frac{\theta_{N_t-1}}{\Delta t^2} \right )
+ \int \zeta \bar E_{N_t-1} (\bar B \bar B^T)^{-1}  \left ( \frac{\theta_{N_t}}{\Delta t^2}
\right ) \nonumber\\
& - \int \zeta (-2 \bar E_1 + \bar E_0) (\bar B \bar B^T)^{-1}  \left ( \frac{\theta_1}{\Delta t^2} \right ) - \int \zeta \bar E_1 (\bar B \bar B^T)^{-1}  \left ( \frac{\theta_2}{\Delta t^2}
\right ) \nonumber.
\end{align}
The boundary terms $\mathcal{B}_1$ and $\mathcal{B}_2$ can be reformulated in the following way. First,
\begin{align*}
\Delta t \mathcal{B}_1
& = 
-\int \Phi(x_{N_t}) \hat \zeta \left ( \sum_{k=3}^{N_t} \bar E_{k-1} A_2^T (\bar B \bar B^T)^{-1} \theta_k \right ) \\
& \qquad\qquad  +
\int \Phi(x_{N_t}) \hat \zeta \left (\sum_{k=2}^{N_t-1} \bar E_{k-1} A_2^T (\bar B \bar B^T)^{-1} \theta_k \right )\\
& = 
-\int \Phi(x_{N_t-1}) \hat \zeta \left ( \sum_{k=2}^{N_t-1} \bar E_{k-1} A_2^T (\bar B \bar B^T)^{-1} \theta_k \right ) \\
& \qquad\qquad  +
\int \Phi(x_{N_t}) \hat \zeta \left (\sum_{k=2}^{N_t-1} \bar E_{k-1} A_2^T (\bar B \bar B^T)^{-1} \theta_k \right ).
\end{align*}
Here we use the stationarity property
and the notation 
$$
\hat \zeta 
= \left ( \prod_{j=1}^{N_t} P_j \right ) P^{ss}.
$$ 
Therefore 
\begin{equation}
\label{eq:B1result}
\mathcal{B}_1 
= \int \hat \zeta \left ( \frac{\Phi(x_{N_t}) - \Phi(x_{N_t-1})}{\Delta t} \right )  
\sum_{k=2}^{N_t-1} \bar E_{k-1} A_2^T (\bar B \bar B^T)^{-1} \theta_k.
\end{equation}
Next, we focus on $\mathcal{B}_2$. By adding and subtracting $\bar E_{N_t}$ and $\bar E_2$ in the second and fourth integrals in $\mathcal{B}_2$, we obtain
\begin{align*}
\mathcal{B}_2
= & \int \zeta (\bar E_{N_t-2} - 2 \bar E_{N_t-1}) (\bar B \bar B^T)^{-1}  \left ( \frac{\theta_{N_t-1}}{\Delta t^2} \right )
+ \int \zeta (\bar E_{N_t-1} - \bar E_{N_t} + \bar E_{N_t})(\bar B \bar B^T)^{-1}  \left ( \frac{\theta_{N_t}}{\Delta t^2}
\right ) \\
& - \int \zeta (-2 \bar E_1 + \bar E_0) (\bar B \bar B^T)^{-1}  \left ( \frac{\theta_1}{\Delta t^2} \right ) 
- \int \zeta (\bar E_1- \bar E_2 + \bar E_2) (\bar B \bar B^T)^{-1}  \left ( \frac{\theta_2}{\Delta t^2}
\right ).
\end{align*}
Then we split into two parts $\mathcal{B}_2 = \mathcal{B}_{2,I} + \mathcal{B}_{2,II}$ where 
\begin{align*}
\mathcal{B}_{2,I}
= & -\int \zeta (\bar E_{N_t-1} - \bar E_{N_t-2}) (\bar B \bar B^T)^{-1}  \left ( \frac{\theta_{N_t-1}}{\Delta t^2} \right )
- \int \zeta (\bar E_{N_t} - \bar E_{N_t-1})(\bar B \bar B^T)^{-1}  \left ( \frac{\theta_{N_t}}{\Delta t^2}
\right ) \\
& + \int \zeta (\bar E_1 - \bar E_0) (\bar B \bar B^T)^{-1}  \left ( \frac{\theta_1}{\Delta t^2} \right ) 
+ \int \zeta (\bar E_2- \bar E_1) (\bar B \bar B^T)^{-1}  \left ( \frac{\theta_2}{\Delta t^2}
\right ),
\end{align*}
and 
\begin{align*}
\mathcal{B}_{2,II}
= &
-\int 
\zeta 
\bar E_{N_t-1}
(\bar B \bar B^T)^{-1}
\left ( \frac{\theta_{N_t-1}}{\Delta t^2} \right )
+\int 
\zeta 
\bar E_{N_t}
(\bar B \bar B^T)^{-1}
\left ( \frac{\theta_{N_t}}{\Delta t^2} \right )\\
& +\int 
\zeta 
\bar E_1
(\bar B \bar B^T)^{-1}
\left ( \frac{\theta_1}{\Delta t^2} \right ) -\int 
\zeta 
\bar E_2
(\bar B \bar B^T)^{-1}
\left ( \frac{\theta_2}{\Delta t^2} \right ).\\
\end{align*}
Next, we reformulate $\mathcal{B}_{2,I}$ and $\mathcal{B}_{2,II}$.
For convenience, we use the notation
$f_k = (\bar E_k - \bar E_{k-1}) (\bar B \bar B^T)^{-1}
\theta_k \Delta t^{-2}$.
We recognize a telescopic sum structure,
\begin{align*}
-\mathcal{B}_{2,I} = &
\int \Phi(x_{N_t}) \hat \zeta \left ( \sum_{k=2}^{N_t-1} f_k - f_{k-1} \right ) 
+ \int \Phi(x_{N_t}) \hat \zeta \left ( \sum_{k=3}^{N_t} f_k - f_{k-1} \right )\\
= &
\int \Phi(x_{N_t}) \hat \zeta \left ( \sum_{k=2}^{N_t-1} f_k \right ) 
- \int \Phi(x_{N_t}) \hat \zeta \left ( \sum_{k=2}^{N_t-1} f_{k-1} \right )\\
& + 
\int \Phi(x_{N_t}) \hat \zeta \left ( \sum_{k=3}^{N_t} f_k \right ) 
- \int \Phi(x_{N_t}) \hat \zeta \left ( \sum_{k=3}^{N_t} f_{k-1} \right ) .
\end{align*}
We use the the stationarity property on the left and a simple change of summation in the right.
\begin{align*}
-\mathcal{B}_{2,I} =
&
\int \Phi(x_{N_t-1}) \hat \zeta \left ( \sum_{k=1}^{N_t-2} f_k \right ) 
- \int \Phi(x_{N_t}) \hat \zeta \left ( \sum_{k=1}^{N_t-2} f_k \right )\\
& + 
\int \Phi(x_{N_t-1}) \hat \zeta \left ( \sum_{k=2}^{N_t-1} f_k \right ) 
- \int \Phi(x_{N_t}) \hat \zeta \left ( \sum_{k=2}^{N_t-1} f_k \right ).
\end{align*}
Therefore we obtain
\begin{align}
\mathcal{B}_{2,I} = & \int \hat \zeta \left ( \frac{\Phi(x_{N_t}) - \Phi(x_{N_t-1})}{\Delta t} \right )  \left ( \sum_{k=1}^{N_t-2} 
\left (\frac{\bar E_k - \bar E_{k-1}}{\Delta t} \right ) (\bar B \bar B^T)^{-1}
\theta_k  \right )
\label{eq:B21result}\\
& + \int \hat \zeta \left ( \frac{\Phi(x_{N_t}) - \Phi(x_{N_t-1})}{\Delta t} \right )  \left ( \sum_{k=2}^{N_t-1} 
\left (\frac{\bar E_k - \bar E_{k-1}}{\Delta t} \right ) (\bar B \bar B^T)^{-1}
\theta_k \right ). \nonumber
\end{align}
Now we turn to the boundary term $\mathcal{B}_{2,II}$ and we use the notation $\mfe_k = \bar E_k
(\bar B \bar B^T)^{-1}
\theta_k \Delta t^{-2}$:
\begin{align*}
\mathcal{B}_{2,II} = &
- \int \zeta \sum_{k=2}^{N_t-1} \left ( \mfe_k - \mfe_{k-1} \right )
+ \int \zeta \sum_{k=3}^{N_t} \left ( \mfe_k - \mfe_{k-1} \right )\\
= &
- \int \Phi(x_{N_t}) \hat \zeta \left ( \sum_{k=2}^{N_t-1} \mfe_k \right ) + \int \Phi(x_{N_t}) \hat \zeta \left ( \sum_{k=2}^{N_t-1} \mfe_{k-1} \right )\\ & + \int \Phi(x_{N_t}) \hat \zeta \left ( \sum_{k=3}^{N_t} \mfe_k \right ) - \int \Phi(x_{N_t}) \hat \zeta \left ( \sum_{k=3}^{N_t} \mfe_{k-1} \right ).
\end{align*}
We use the stationary property and elementary changes of summation to get
\begin{align*}
\mathcal{B}_{2,II} = &
- \int \Phi(x_{N_t-1}) \hat \zeta \left ( \sum_{k=1}^{N_t-2} \mfe_k \right ) + \int \Phi(x_{N_t}) \hat \zeta \left ( \sum_{k=1}^{N_t-2} \mfe_k \right )\\ & + \int \Phi(x_{N_t-2}) \hat \zeta \left ( \sum_{k=1}^{N_t-2} \mfe_k \right ) - \int \Phi(x_{N_t-1}) \hat \zeta \left ( \sum_{k=1}^{N_t-2} \mfe_{k} \right ).\\
\end{align*}
Therefore,
\begin{equation}
\label{eq:B22result}
\mathcal{B}_{2,II}
=
\int \hat \zeta \left (\frac{\Phi(x_{N_t}) - 2 \Phi(x_{N_t-1}) + \Phi(x_{N_t-2})}{\Delta t^2} 
\right )  \left ( \sum_{k=1}^{N_t-2} \bar E_k
(\bar B \bar B^T)^{-1}
\theta_k \right ).
\end{equation}
Now we collect all our results for the boundary terms  
\eqref{eq:B1result},
\eqref{eq:B21result},
\eqref{eq:B22result}
and substitute them into the right hand side of 
\eqref{eq:sumTi}
before passing to the limit.
We obtain
\begin{align}\label{eq:sumTi_reform1}
\sum_{i=3}^{N_t} \mathcal{T}_i
= & \sum_{i=3}^{N_t} \int \hat \zeta \Phi(x_{N_t}) \bar E_{i-1}
A_1^T (\bar B \bar B^T)^{-1}
\theta_i\\
& + \sum_{i=2}^{N_t-1} \int \hat \zeta \Phi(x_{N_t}) \left ( \frac{\bar E_i - \bar E_{i-1}}{\Delta t} \right ) A_2^T (\bar B \bar B^T)^{-1}
\theta_i \nonumber\\
& + \sum_{i=1}^{N_t-2} \int \hat \zeta \Phi(x_{N_t}) \left ( \frac{\bar E_{i+1}- 2 \bar E_i + \bar E_{i-1}}{\Delta t^2}  \right )
(\bar B \bar B^T)^{-1} \theta_i \nonumber\\
& + \int \hat \zeta \left ( \frac{\Phi(x_{N_t}) - \Phi(x_{N_t-1})}{\Delta t} \right )  
\sum_{k=2}^{N_t-1} \bar E_{k-1} A_2^T (\bar B \bar B^T)^{-1} \theta_k \nonumber\\
& + \int \hat \zeta \left ( \frac{\Phi(x_{N_t}) - \Phi(x_{N_t-1})}{\Delta t} \right ) \left ( \sum_{k=1}^{N_t-2} 
\left (\frac{\bar E_k - \bar E_{k-1}}{\Delta t} \right ) (\bar B \bar B^T)^{-1}
\theta_k \right )
\nonumber \\
& + \int \hat \zeta \left ( \frac{\Phi(x_{N_t}) - \Phi(x_{N_t-1})}{\Delta t} \right ) \left ( \sum_{k=2}^{N_t-1} 
\left (\frac{\bar E_k - \bar E_{k-1}}{\Delta t} \right ) (\bar B \bar B^T)^{-1}
\theta_k \right )\nonumber \\
& + \int \hat \zeta \left (\frac{\Phi(x_{N_t}) - 2 \Phi(x_{N_t-1}) + \Phi(x_{N_t-2})}{\Delta t^2} 
\right ) \left ( \sum_{k=1}^{N_t-2} \bar E_k
(\bar B \bar B^T)^{-1}
\theta_k \right ). \nonumber
\end{align}
Again assuming we can take the limit as $N_t \to \infty$, we obtain our second main result when $n'=2$, i.e.
\begin{align*}
     \frac{\dd}{\dd\lambda} \langle \Phi(x(t)) \rangle_\lambda \Big|_{\lambda=0}= & \:
      \Big\langle \Phi(x(t)) \big[ p_{00}(t) + p_{10}(t) + p_{20}(t) \big] \Big\rangle
    + \Big\langle \dot \Phi(x(t)) \big[ p_{11}(t) + 2 p_{21}(t) \big] \Big\rangle\\
     & + \Big\langle \ddot \Phi(x(t)) \; p_{22}(t) \Big\rangle,
\end{align*}
where
$$
\dot p_{0,0}(t) = \bar E(x(t)) A_1^T (\bar B \bar B^T)^{-1} B \dot w(t), \quad
p_{1,0}(t) = \frac{\dd}{ \dd t} \left ( \bar E(x(t)) \right ) A_2^T (\bar B \bar B^T)^{-1} B \dot w(t), 
$$
$$
\dot p_{2,0}(t) = \frac{\dd^2}{ \dd t^2} \left ( \bar E(x(t)) \right ) (\bar B \bar B^T)^{-1} B \dot w(t),
\quad
\dot p_{1,1}(t) =
\bar E(x(t)) A_2^T (\bar B \bar B^T)^{-1} B \dot w(t),
$$
$$ 
\dot p_{2,1}(t) =
\frac{\dd}{ \dd t} \left ( \bar E(x(t)) \right ) 
(\bar B \bar B^T)^{-1} B \dot w(t), \: \mbox{ and } \:
\dot p_{2,2}(t) =
\bar E(x(t)) (\bar B \bar B^T)^{-1} B \dot w(t).
$$

\subsubsection{Case $n'\geq 1$ : key steps }
Now we turn to the general case $n'\geq 1$.\\

\noindent 
\textbf{Step 3}\\
The counterpart 
of formula \eqref{eq:Ti_reform2}
becomes
\begin{equation}
\label{eq:Ti_expansion_gen}
    \sum \limits_{i=n'+1}^{N_t} \mathcal{T}_i
    = \sum \limits_{k=0}^{n'} \underbrace{\frac{(-1)^{k}}{( \Delta t )^{k}} 
    \left ( \sum \limits_{i=n'+1}^{N_t} \int  
    \bar E(x_{i-1}) A_{k+1}^T (\bar B \bar B^T)^{-1}
    (D_k \theta)_i
    \Phi(x_{N_t}) 
    \hat \zeta \right ).}_{S_k :=}
\end{equation}
The finite difference at the order $k \geq 0$ is defined as follows:
$(D_0 \theta)_i = \theta_i$ and $(D_k \theta)_i = (D_{k-1} \theta)_i
- (D_{k-1} \theta)_{i-1}$ for $k \geq 1$.
We recall that the matrices  
$A_k$ for $1 \leq k \leq n'$ are defined in \eqref{eq:ABthem2} and $A_{n'+1} = I_{q'}$.\\

Next we use a summation by part formula shown in Lemma \ref{appendix:lemma_sbp} to reformulate each of the $S_k$, that is swapping the finite difference on $\theta$
to a finite difference on $\bar E$ and obtaining boundary terms. \\

\noindent \textbf{Step 4}\\
Below, we use the notation
$
h_{i,k} =
A_k^T (\bar B \bar B^T)^{-1} \theta_i.
$
In the rhs of  \eqref{eq:Ti_expansion_gen}, we consider each of the $n'+1$ sums individually.
For $k = 0$, we can rewrite the sum as 
$$
S_0 = \sum \limits_{i=n'+1}^{N_t} \int \hat \zeta \Phi(x_{N_t}) \bar E_{i-1} h_{i,1}
$$
and for $k \geq 1$, we apply Lemma \ref{appendix:lemma_sbp} (summation by part in appendix) to obtain $S_k = M_0 + B_0$
where
\begin{align*}
 M_0=  & \sum_{i=n'+1}^{N_t} \int \hat \zeta \Phi(x_{N_t}) \bigg( \frac{D_k \bar E_{i-1}}{(\Delta t)^{k}} \bigg) h_{i-k,k+1},\\ 
B_0 = & \frac{1}{(\Delta t)^{k}} \sum_{j_1=1}^k (-1)^{k+j_1} \int \zeta \bigg( D_{j_1-1} \bar E_{i-1} \bigg) \bigg( D_{k-j_1} h_{i-j_1+1,k+1} \bigg) \bigg |_{n'}^{N_t}. 
\end{align*}
Then, we use the stationarity property to transform $B_0$ into
\begin{align*}
    B_0 = & \frac{1}{(\Delta t)^{k}} \sum_{j_1=1}^k (-1)^{k+j_1} \sum_{i=n'}^{{N_t}-1} \int \hat \zeta \bigg( D_1 \Phi_{N_t} \bigg) \bigg( D_{j_1-1} \bar E_{i-1} \bigg) \bigg( D_{k-j_1} h_{i-j_1+1,k+1} \bigg).
\end{align*}
\textbf{Step 5}\\
We decompose $B_0 = M_{1,1}+M_{1,2}+B_1$ by separating the sum into two parts ($k-j_1=0$ and $k-j_1 \neq 0$):
\begin{align*}
    M_{1,1} = &  \sum_{\substack{j_1=1, \\ j_1 = k}}^{k} \sum_{i=n'}^{N_t-1} \int \hat \zeta \bigg( \frac{D_1 \Phi_{N_t}}{(\Delta t)^{1}} \bigg) \bigg( \frac{D_{k-1} \bar E_{i-1}}{(\Delta t)^{k-1}} \bigg) h_{i-k+1,k+1},\\
    M_{1,2} = &  \sum_{\substack{j_1=1, \\ j_1 \neq k}}^{k} \sum_{i=n'}^{{N_t}-1} \int \hat \zeta \bigg( \frac{D_1 \Phi_{N_t}}{(\Delta t)^{1}} \bigg) \bigg( \frac{D_{k-1} \bar E_{i-1}}{(\Delta t)^{k-1}} \bigg) h_{i-k+1,k+1} ,
    \\ 
    B_1 = &  \frac{1}{(\Delta t)^{k}} \sum_{j_1=1}^{k-1} \sum_{j_2=1}^{k-j_1} (-1)^{k+j_1+j_2} \\
    & \times \int \hat \zeta \bigg( D_1 \Phi_{N_t} \bigg) \bigg( D_{j_1-j_2-2} \bar E_{i-1} \bigg) \bigg( D_{k-j_1-j_2} h_{i-j_1-j_2+2,k+1} \bigg) \bigg|_{n'-1}^{{N_t}-1}. 
\end{align*}
Notice that 
\begin{align*}
    M_{1,1} + M_{1,2} = 
    & \sum_{j_1=1}^{k} \sum_{i=n'}^{N_t-1} \int \hat \zeta \bigg( \frac{D_1 \Phi_{N_t}}{(\Delta t)^{1}} \bigg) \bigg( \frac{D_{k-1} \bar E_{i-1}}{(\Delta t)^{k-1}} \bigg) h_{i-k+1,k+1}.
\end{align*}

\noindent \textbf{Step 6}\\
We repeat the steps above to the boundary term $(B_1)$ and the subsequent boundary terms that appear at each iteration.
The result obtained after repeating these a total of $k-1$ times is
\begin{align*}
    & \sum_{m=1}^{k-1} \binom{k}{m} \sum_{i=n'-(m-1)}^{N_t-m} \int \hat \zeta \bigg( \frac{D_m \Phi_{N_t}}{(\Delta t)^k} \bigg) \bigg( \frac{D_{k-m} \bar E_{i-1}}{(\Delta t)^k} \bigg) h_{i-k+m,k+1}
    \\& + \sum_{j_1=1}^{1} \sum_{j_2=1}^{2-j_1} ... \sum_{j_{k}=1}^{k-\mathcal{J}_{m-1}} (-1)^{k+\mathcal{J}_m} \int \hat \zeta \bigg( D_{k-1} \Phi_{N_t} \bigg) \bigg( D_{\mathcal{J}_k-k} \bar E_{i-1} \bigg) \bigg( D_{k-\mathcal{J}_k} h_{i-\mathcal{J}_k+k,k+1} \bigg) \bigg|_{n'-(l-1)}^{N_t-(l-1)}.
\end{align*}
Here, $\mathcal{J}_k=j_1+j_2+...+j_k$.\\

\noindent \textbf{Step 7}\\
By applying the stationarity property to the boundary term, collecting terms, and rearranging, we arrive at the desired result
\begin{equation}
\sum_{i=n'+1}^{N_t} \mathcal{T}_i
 = \sum_{k=0}^{n'} \sum_{m=0}^k \binom{k}{m} \sum_{i=n'+1-k}^{N_t-k} \int \hat \zeta \bigg( \frac{D_m \Phi_{N_t}}{(\Delta t)^m} \bigg) \bigg( \frac{D_{k-m} \bar E_{i-1+(k-m)}}{(\Delta t)^{k-m}} \bigg) A_{k+1}^T (\bar B \bar B^T)^{-1} \theta_{i}.
\end{equation}

\noindent \textbf{Step 8}\\
Taking limits as $N_t \to \infty$ and rearranging further, we then have the formula 
(\ref{eq:sensitivity2}-\ref{eq:weight4formula2}).
\subsubsection{ Technical lemmata}
\begin{lemma}
\label{appendix:lemma_3132}
$\mathcal{T}_i$ defined by \eqref{eq:Ti}
satisfies 
\eqref{eq:Ti_reform1}.
\end{lemma}
\begin{proof}
We start from \eqref{eq:Ti} with $n'=2$,
\begin{align*}
\mathcal{T}_i & = \int \Phi(x_{N_t}) P_{\lambda=0,i}' \left ( \prod_{j=1, j \neq i}^{N_t} P_j \right ) P^{ss}\\
& = \int \Phi(x_{N_t}) 
\bar E_{i-1} \cdot
\left ( \delta_{\mathcal{X}_i}\right )_{y_{i-1,1}} 
\left ( \delta_{\mathcal{Y}_{i-1,1}}  \delta_{\mathcal{Y}_{i,1}} 
g_i \right)
\left ( \prod_{j=1, j \notin \{i-1,i\}}^{N_t} P_j \right ) \delta_{\mathcal{X}_{i-1}} g_{i-1} P^{ss}.
\end{align*}
An integration by part gives 
$\mathcal{T}_i =  
- \mathcal{T}_{i,1} 
- \mathcal{T}_{i,2}
- \mathcal{T}_{i,3}$
where 
\begin{align*}
\mathcal{T}_{i,1} & = 
\int \Phi(x_{N_t}) 
\bar E_{i-1} \cdot 
\left ( \delta_{\mathcal{Y}_{i-1,1}} \right )_{y_{i-1,1}} 
\left ( \prod_{j=1, j \neq i-1}^{N_t} P_j \right ) \delta_{\mathcal{X}_{i-1}} g_{i-1} P^{ss}\\
\mathcal{T}_{i,2} & = 
\int \Phi(x_{N_t}) 
\bar E_{i-1} \cdot 
\left ( \delta_{\mathcal{Y}_{i,1}} \right )_{y_{i-1,1}}
\delta_{\mathcal{X}_i} g_i
\left ( \prod_{j=1, j \neq i}^{N_t} P_j \right ) P^{ss}\\
\mathcal{T}_{i,3} & = 
\int \Phi(x_{N_t}) 
\bar E_{i-1} \cdot 
\left ( g_i \right )_{y_{i-1,1}}
\delta_{\mathcal{Y}_{i,1}}
\delta_{\mathcal{X}_i} 
\left ( \prod_{j=1, j \neq i}^{N_t} P_j \right ) P^{ss}.
\end{align*}
Using 
$ 
\left ( \delta_{\mathcal{Y}_{i-1,1}} \right )_{y_{i-1,1}} 
= - (\Delta t)^{-1}
\left ( \delta_{\mathcal{Y}_{i-1,1}} \right )_{y_{i-2,2}}
$, we reformulate $\mathcal{T}_{i,1}$ as follows
\begin{align*}
\mathcal{T}_{i,1} 
& = - \frac{1}{\Delta t}
\int \Phi(x_{N_t}) 
\bar E_{i-1} \cdot
\left ( \delta_{\mathcal{Y}_{i-1,1}} \right )_{y_{i-2,2}} 
\left ( \prod_{j=1, j \neq i-1}^{N_t} P_j \right ) \delta_{\mathcal{X}_{i-1}} g_{i-1} P^{ss}\\
& = 
\frac{1}{\Delta t}
\int \zeta 
\bar E_{i-1} \cdot 
\left ( \log( g_{i-2} g_{i-1} ) \right)_{y_{i-2,2}} 
\: \mbox{(by integration by parts)},
\end{align*}
where we reuse the notation
$$
\zeta 
= \Phi(x_{N_t})
\left ( \prod_{j=1}^{N_t} P_j \right ) P^{ss}.
$$
Next, using 
$ 
\left ( \delta_{\mathcal{Y}_{i,1}} \right )_{y_{i-1,1}} 
= (\Delta t)^{-1} 
\left ( \delta_{\mathcal{Y}_{i-1,1}} \right )_{y_{i-1,2}}
$, we reformulate $\mathcal{T}_{i,2}$ as follows
\begin{align*}
\mathcal{T}_{i,2} & = 
\frac{1}{\Delta t}
\int \Phi(x_{N_t}) 
\bar E_{i-1} \cdot
 \left ( \delta_{\mathcal{Y}_{i,1}} \right )_{y_{i-1,2}} 
\left ( \prod_{j=1, j \neq i }^{N_t} P_j \right ) \delta_{\mathcal{X}_{i}} g_{i} P^{ss}\\
& = 
-\frac{1}{\Delta t}
\int \zeta 
\bar E_{i-1} \cdot 
\left ( \log( g_{i-1} g_i ) \right)_{y_{i-1,2}} 
\: \mbox{(by integration by parts)}.
\end{align*}
Finally, we reformulate $\mathcal{T}_{i,3}$ as follows
$$
\mathcal{T}_{i,3} 
=
\int \zeta 
\bar E_{i-1} \cdot 
\left ( \log (g_i) \right )_{y_{i-1,1}}.
$$
We can collect the results above to obtain \eqref{eq:Ti_reform1}.
\end{proof}
\begin{lemma}
\label{appendix:lemma_sbp}
Summation by parts formula for $1 \leq k \leq n$
$$
    \sum_{i=n+1}^N \phi_{i-1} D_k \psi_i = (-1)^k \sum_{i=n+1}^N D_k \phi_{i-1} \psi_{i-k} + \sum_{j=1}^k (-1)^{j-1} \bigg( D_{j-1} \phi_{i-1} D_{k-j} \psi_{i-j+1} \bigg) \bigg |_n^N.
$$
\end{lemma}
\begin{proof}
We have that
\begin{align*}
    \sum_{i=n+1}^N \phi_{i-1} D^k \psi_i
    & = - \bigg( \sum_{i=n+1}^N \phi_{i-1} D^{k-1} \psi_{i-1} - \sum_{i=n+1}^N \phi_{i-2} D^{k-1} \psi_{i-1} \bigg) \\
    & \qquad + \bigg( \sum_{i=n+1}^N \phi_{i-1} D^{k-1} \psi_i - \sum_{i=n+1}^N \phi_{i-2} D^{k-1} \psi_{i-1} \bigg)
    \\& = - \sum_{i=n+1}^N D^1 \phi_{i-1} D^{k-1} \psi_{i-1} + \bigg( \phi_{i-1} D^{k-1} \psi_i \bigg) \bigg|_n^N.
\end{align*}
By inductively applying the above identity to $\sum_{i=n+1}^N D^l \phi_{i-1} D^{k-l} \psi_{i-l}$, it can be shown that the following holds for $l=1,...,k$.
\begin{align*}
    \sum_{i=n+1}^N \phi_{i-1} D^k \psi_i = (-1)^l \sum_{i=n+1}^N D^l \phi_{i-1} D^{k-l} \psi_{i-l} + \sum_{j=1}^l (-1)^{j-1} \bigg( D^{j-1} \phi_{i-1} D^{k-j} \psi_{i-j+1} \bigg) \bigg |_n^N
\end{align*}
In particular, the case $l=k$ holds and concludes the proof.
\end{proof}
\subsection{Details of weight formulas}
\label{appendix:weights}
Calculations yield the following formulae. In the rest of this section, let $\alpha = k/\xi_0 + \ell$ and $\beta = k/\xi_0 - \ell$. Furthermore, when the average, $\avg{...}$, appears without the dependence on $t$, we implicitly mean to take the average with respect to the steady state.
\subsubsection{Weight formulas associated with \eqref{eq:decomp1}}
$$
    \avg{ x(t) p_{0,0}^1(t) } = \frac{1 - e^{-\frac{k}{\xi_0} t}}{k}  - \frac{e \sigma}{\xi_0} c(t),
    \quad \mbox{ and } \quad 
    \avg{ \dot x(t) p_{1,1}^1(t) } = \frac{e \sigma}{\xi_0} c(t),
$$
where
\begin{fleqn}
$$
c(t) = \frac{e^{-\ell t}(\beta\cos(\omega t) + \omega \sin(\omega t)) - \beta e^{-\frac{k}{\xi_0} t}}{\beta^2+\omega^2}.
$$
\end{fleqn}

\subsubsection{Weight formulas associated with \eqref{eq:decomp2}}
\begin{align*}
    \avg{x^2(t)p_{0,0}^2(t)} & = \frac{\sigma^2 }{\xi_0^2} e^{-2\ell t} \bigg( \frac{\cos(2\omega t) \mathcal{A}_1 + \sin(2\omega t) \mathcal{B}_1 + \mathcal{C}_1}{2\beta(\beta^2+\omega^2)^2} \bigg)
    \\& + \frac{2\sigma}{\xi_0} e^{-\alpha t} \bigg( \frac{\cos(\omega t)\mathcal{D}_1 + \sin(\omega t)\mathcal{E}_1}{\beta^2+\omega^2} \bigg)
    \\& + \frac{\sigma}{k} (e^{-2k/\xi_0 t} - 1) a_{1,1},
\end{align*}
\begin{align*}
    \avg{x^2(t)p_{1,0}^2(t)} & = \frac{\sigma^2 }{\xi_0^2} e^{-2\ell t} \bigg( \frac{\cos(2\omega t) \mathcal{A}_2 + \sin(2\omega t) \mathcal{B}_2 + \mathcal{C}_2}{2(\frac{k}{\xi_0}-\ell)((\frac{k}{\xi_0}-\ell)^2+\omega^2)^2} \bigg)
    \\& + \frac{2\sigma}{\xi_0} e^{-(k/\xi_0+\ell)t} \bigg( \frac{\cos(\omega t)\mathcal{D}_2 + \sin(\omega t)\mathcal{E}_2}{(\frac{k}{\xi_0}-\ell)^2+\omega^2} \bigg)
    \\& + \frac{\sigma}{k} (e^{-2k/\xi_0 t} - 1) a_{1,2},
\end{align*}
and
\begin{align*}
    \avg{\dot{x^2}(t)p_{1,1}^2(t)} & = \frac{-2k\sigma^2 }{\xi_0^3} e^{-2\ell t} \bigg( \frac{\cos(2\omega t) \mathcal{A}_3 + \sin(2\omega t) \mathcal{B}_3 + \mathcal{C}_3}{2(\frac{k}{\xi_0}-\ell)((\frac{k}{\xi_0}-\ell)^2+\omega^2)^2} \bigg)
    \\& + \frac{-4k\sigma}{\xi_0^2} e^{-(k/\xi_0+\ell)t} \bigg( \frac{\cos(\omega t)\mathcal{D}_3 + \sin(\omega t)\mathcal{E}_3}{(\frac{k}{\xi_0}-\ell)^2+\omega^2} \bigg)
    \\& + \frac{-2\sigma}{\xi_0} (e^{-2k/\xi_0 t} - 1) a_{1,3}
    \\& + \frac{\sigma}{\xi_0} e^{-2\ell t} \bigg(\frac{\cos^2(\omega t)\mathcal{F} + \cos(\omega t)\sin(\omega t)\mathcal{G} + \sin^2(\omega t)\mathcal{H}}{\alpha^2+\omega^2} \bigg)
    \\& + e^{-(k/\xi_0+\ell)t} (\cos(\omega t)a_{1,3} + \sin(\omega t)a_{2,3}) - a_{1,3},
\end{align*}
where we use the following constant
\begin{fleqn}
\begin{align*}
    & \mathcal{A}_i = c_{1,i}\beta^3 - c_{2,i}\beta\omega^2 - (2b_{2,i})\beta^2\omega,\\
    & \mathcal{B}_i = c_{1,i}\beta^2\omega + c_{2,i}\beta^2\omega + (2b_{2,i})\beta^3,\\
    & \mathcal{C}_i = c_{2,i}\beta(\beta^2+\omega^2) + (2b_{2,i})\omega(\beta^2+\omega^2),\\
    & \mathcal{D}_i = a_{1,i}\beta - a_{2,i}\omega,\\
    & \mathcal{E}_i = a_{1,i}\omega+a_{2,i}\beta, \\
    & \mathcal{F} = (b_{1,3}a+b_{2,3}\omega),\\
    & \mathcal{G} = (b_{1,3}\omega+2b_{2,3}a+b_{3,3}\omega),
    \\
    & \mathcal{H} = (b_{2,3}\omega+b_{3,3}a),
\end{align*}
and
\begin{align*}
    c_{1,i} &= \frac{4\ell^2\avg{y_1F_{1,i}} + 4\omega\ell\avg{y_2F_{1,i}} + 4\omega\ell\avg{y_1F_{2,i}} - 4\ell^2\avg{y_2F_{2,i}}}{-4(\ell^2+\omega^2)},
    \\ c_{2,i} &= -(\avg{y_1F_{1,i}}+\avg{y_2F_{2,i}}),
    \\ b_{1,i} &= \frac{2(2\ell^2+\omega^2)\avg{y_1F_{1,i}} + 2\omega\ell\avg{y_2F_{1,i}} + 2\omega\ell\avg{y_1F_{2,i}} + 2\omega^2\avg{y_2F_{2,i}}}{-4(\ell^2+\omega^2)},
    \\ b_{2,i} &= \frac{-2\omega\ell\avg{y_1F_{1,i}} + 2\ell^2\avg{y_2F_{1,i}} + 2\ell^2\avg{y_1F_{2,i}} + 2\omega\ell\avg{y_2F_{2,i}}}{-4(\ell^2+\omega^2)},
    \\ b_3 &= \frac{\omega^22\avg{y_1F_{1,i}} + -2\ell\omega(\avg{y_2F_{1,i}} + \avg{y_1F_{2,i}}) + (2\ell^2+\omega^2)2\avg{y_2F_{2,i}}}{-4(\ell^2+\omega^2)},
    \\ a_{1,i} &= \frac{-\alpha(\ell\avg{xF_{1,i}}-b_{1,i}) + \omega(\ell\avg{xF_{2,i}}-b_{2,i})},{\alpha^2+\omega^2},
    \\ a_{2,i} &= \frac{\omega(\ell\avg{xF_{1,i}}-b_{1,i}) - \alpha(\ell\avg{xF_{2,i}}-b_{2,i})}{\alpha^2+\omega^2},
\end{align*}
\end{fleqn}
with
$
F_{1,1} = e x, \:
F_{2,1} = e \omega \ell^{-1} x, \:
F_{1,2} = e\ell^{-1} \dot x, \:
F_{2,2} = 0, \:
F_{1,3} = e \ell^{-1} x, \:
F_{2,3} = 0.
$

\subsubsection{Weight formulas associated with \eqref{eq:decomp3}}
\begin{align*}
    \avg{x(t)p_{00}^3(t)} & = e\gamma \alpha_0 \bigg( \frac{1 - e^{-k/\xi_0 t}}{k} - \ell \frac{te^{-\ell t}}{k-\ell\xi_0} - (k-2\ell\xi_0) \frac{e^{-\ell t} - e^{-k/\xi_0 t}}{(k-\ell\xi_0)^2} \bigg),
\end{align*}
\begin{align*}
    \avg{\dot x(t) p_{11}^3(t)} & = e\gamma \alpha_1 \bigg( \ell \frac{te^{-\ell t}}{k-\ell\xi_0}
    - \ell \xi_0 \frac{e^{-\ell t} - e^{-k/\xi_0 t}}{(k-\ell\xi_0)^2} \bigg),
\end{align*}
and 
\begin{align*}
    \avg{\ddot x(t) p_{22}^3(t)} & = e \gamma \bigg( -\ell\frac{t e^{-\ell t}}{k-\ell\xi_0} 
    + k\frac{e^{-\ell t} - e^{-k/\xi_0 t}}{(k-\ell\xi_0)^2} \bigg).
\end{align*}

\subsubsection{Weight formulas associated with \eqref{eq:decomp4}}
\begin{fleqn}
\begin{align*}
    \avg{x^2(t)p_{0,0}^4(t)} & = \mathcal{A}_1 \bigg( t^2 e^{-2\ell t} \bigg)
    + \mathcal{B}_1 \bigg( t e^{-2\ell t} \bigg)
    + \mathcal{C}_1 \bigg( t e^{-\ell t} \bigg)
    + \mathcal{D}_1 \bigg( t e^{-\alpha t} \bigg)
    \\& + \mathcal{E}_1 \bigg( e^{-2\ell t} - e^{-2k/\xi_0t} \bigg)
    + \mathcal{F}_1 \bigg( e^{-\ell t} - e^{-2k/\xi_0t} \bigg)
    + \mathcal{G}_1 \bigg( e^{-\alpha t} - e^{-2k/\xi_0t} \bigg)
    \\& + \frac{\gamma}{k} b_{1,1} \bigg( 1 - e^{-2k/\xi_0t} \bigg),
\end{align*}
\begin{align*}
    \avg{x^2(t)p_{1,0}^4(t)} & = \mathcal{A}_2 \bigg( t^2 e^{-2\ell t} \bigg)
    + \mathcal{B}_2 \bigg( t e^{-2\ell t} \bigg)
    + \mathcal{C}_2 \bigg( t e^{-\ell t} \bigg)
    + \mathcal{D}_2 \bigg( t e^{-\alpha t} \bigg)
    \\& + \mathcal{E}_2 \bigg( e^{-2\ell t} - e^{-2k/\xi_0t} \bigg)
    + \mathcal{F}_2 \bigg( e^{-\ell t} - e^{-2k/\xi_0t} \bigg)
    + \mathcal{G}_2 \bigg( e^{-\alpha t} - e^{-2k/\xi_0t} \bigg)
    \\& + \frac{\gamma}{k} b_{1,2} \bigg( 1 - e^{-2k/\xi_0t} \bigg),
\end{align*}
\begin{align*}
    \avg{x^2(t)p_{1,0}^4(t)} & = \mathcal{A}_3 \bigg( t^2 e^{-2\ell t} \bigg)
    + \mathcal{B}_3 \bigg( t e^{-2\ell t} \bigg)
    + \mathcal{C}_3 \bigg( t e^{-\ell t} \bigg)
    + \mathcal{D}_3 \bigg( t e^{-\alpha t} \bigg)
    \\& + \mathcal{E}_3 \bigg( e^{-2\ell t} - e^{-2k/\xi_0t} \bigg)
    + \mathcal{F}_3 \bigg( e^{-\ell t} - e^{-2k/\xi_0t} \bigg)
    + \mathcal{G}_3 \bigg( e^{-\alpha t} - e^{-2k/\xi_0t} \bigg)
    \\& + \frac{\gamma}{k} b_{1,3} \bigg( 1 - e^{-2k/\xi_0t} \bigg),
\end{align*}
\begin{align*}
    \avg{\dot{x^2}(t)p_{1,1}^4(t)} & = ( \frac{-2k}{\xi_0} \mathcal{A}_4 + \frac{2\gamma}{\xi_0} \mathcal{I}_4) \bigg( t^2 e^{-2\ell t} \bigg)
    + (\frac{-2k}{\xi_0} \mathcal{B}_4 + \frac{2\gamma}{\xi_0}\mathcal{J}_4) \bigg( t e^{-2\ell t} \bigg)
    \\& + \frac{-2k}{\xi_0} \mathcal{C}_4 \bigg( t e^{-\ell t} \bigg)
    + (\frac{-2k}{\xi_0} \mathcal{D}_4 + \frac{2\gamma}{\xi_0}\mathcal{L}_4) \bigg( t e^{-\alpha t} \bigg)
    \\& + \frac{-2k}{\xi_0} \mathcal{E}_4 \bigg( e^{-2\ell t} - e^{-2k/\xi_0t} \bigg)
    + \frac{-2k}{\xi_0} \mathcal{F}_4 \bigg( e^{-\ell t} - e^{-2k/\xi_0t} \bigg)
    \\& + \frac{-2k}{\xi_0} \mathcal{G}_4 \bigg( e^{-\alpha t} - e^{-2k/\xi_0t} \bigg)
    + \frac{2\gamma}{\xi_0}b_{1,4} \bigg( e^{-2k/\xi_0t} \bigg)
    \\& + \frac{2\gamma}{\xi_0}\mathcal{K}_4 \bigg( e^{-2\ell t} \bigg)
    + \frac{2\gamma}{\xi_0}\mathcal{M}_4 \bigg( e^{-\alpha t} \bigg),
\end{align*}
\begin{align*}
    \avg{\dot{x^2}(t)p_{2,1}^4(t)} & = ( \frac{-2k}{\xi_0} \mathcal{A}_5 + \frac{2\gamma}{\xi_0} \mathcal{I}_5) \bigg( t^2 e^{-2\ell t} \bigg)
    + (\frac{-2k}{\xi_0} \mathcal{B}_5 + \frac{2\gamma}{\xi_0}\mathcal{J}_5) \bigg( t e^{-2\ell t} \bigg)
    \\& + \frac{-2k}{\xi_0} \mathcal{C}_5 \bigg( t e^{-\ell t} \bigg)
    + (\frac{-2k}{\xi_0} \mathcal{D}_5 + \frac{2\gamma}{\xi_0}\mathcal{L}_5) \bigg( t e^{-\alpha t} \bigg)
    \\& + \frac{-2k}{\xi_0} \mathcal{E}_5 \bigg( e^{-2\ell t} - e^{-2k/\xi_0t} \bigg)
    + \frac{-2k}{\xi_0} \mathcal{F}_5 \bigg( e^{-\ell t} - e^{-2k/\xi_0t} \bigg)
    \\& + \frac{-2k}{\xi_0} \mathcal{G}_5 \bigg( e^{-\alpha t} - e^{-2k/\xi_0t} \bigg)
    + \frac{2\gamma}{\xi_0}b_{1,5} \bigg( e^{-2k/\xi_0t} \bigg)
    \\& + \frac{2\gamma}{\xi_0}\mathcal{K}_5 \bigg( e^{-2\ell t} \bigg)
    + \frac{2\gamma}{\xi_0}\mathcal{M}_5 \bigg( e^{-\alpha t} \bigg),
\end{align*}
and
\begin{align*}
    \avg{\ddot{x^2}(t)p_{2,2}^4(t)} & = (\frac{4k^2}{\xi_0^2}\mathcal{A}_6 - \frac{6k\sigma}{\xi_0^2}\mathcal{I}_6 + \frac{2\sigma^2}{\xi_0^2}\mathcal{N}_6 - \frac{2\sigma}{\xi_0}\mathcal{I}_6) \bigg( t^2 e^{-2\ell t} \bigg)
    \\& + (\frac{4k^2}{\xi_0^2}\mathcal{B}_6 - \frac{6k\sigma}{\xi_0^2}\mathcal{J}_6 + \frac{2\sigma^2}{\xi_0^2}\mathcal{O}_6 + \frac{2\sigma}{\xi_0}\mathcal{R}_6) \bigg( t e^{-2\ell t} \bigg)
    \\& + \frac{4k^2}{\xi_0^2}\mathcal{C}_6 \bigg( t e^{-\ell t} \bigg)
    + (\frac{4k^2}{\xi_0^2}\mathcal{D}_6 - \frac{6k\sigma}{\xi_0^2}\mathcal{L}_6 - \frac{2\sigma}{\xi_0}\mathcal{L}_6) \bigg( t e^{-\alpha t} \bigg)
    \\& + \frac{4k^2}{\xi_0^2}\mathcal{E}_6 \bigg( e^{-2\ell t} - e^{-2k/\xi_0t} \bigg)
    + \frac{4k^2}{\xi_0^2}\mathcal{F}_6 \bigg( e^{-\ell t} - e^{-2k/\xi_0t} \bigg)
    \\& + \frac{4k^2}{\xi_0^2}\mathcal{G}_6 \bigg( e^{-\alpha t} - e^{-2k/\xi_0t} \bigg)
    + \frac{4k^2}{\xi_0^2}\frac{\gamma}{k} b_{1,i} \bigg( 1 - e^{-2k/\xi_0t} \bigg)
    \\& + (- \frac{6k\sigma}{\xi_0^2}\mathcal{K}_6 - \frac{2\sigma}{\xi_0}\frac{1}{2} \mathcal{K}_6) \bigg( e^{-2\ell t} \bigg)
    + (- \frac{6k\sigma}{\xi_0^2}\mathcal{M}_6 + \frac{2\sigma}{\xi_0}\mathcal{U}_6) \bigg( e^{-\alpha t} \bigg)
    \\& + (- \frac{6k\sigma}{\xi_0^2}b_{1,6} + \frac{2\sigma^2}{\xi_0^2}\mathcal{P}_6 + \frac{2\sigma}{\xi_0}b_{2,6})
\end{align*}
\end{fleqn}
where we use the following constants
\begin{fleqn}
\begin{align*}
    & \mathcal{A}_i = \frac{\sigma}{\xi_0}A_i(-1/\beta^2),
    \\& \mathcal{B}_i = \frac{\sigma}{\xi_0}A_i(1/\beta^3-2/(\ell\beta^2)),
    \\& \mathcal{C}_i = \frac{\sigma}{\xi_0}A_i(1/\beta^3),
    \\& \mathcal{D}_i = (\frac{\sigma}{\xi_0}A_i(-2/\beta^3) + 2\gamma/\xi_0 \ell(b_{1,i}+b_{2,i})/\beta),
    \\& \mathcal{E}_i = \frac{\sigma}{\xi_0}A_i(-1/(2\beta^4)+1/(\ell\beta^3)),
    \\& \mathcal{F}_i = \frac{\sigma}{\xi_0}A_i(-1/(2\beta^4)+2/(\ell\beta^3)),
    \\& \mathcal{G}_i = (\frac{\sigma}{\xi_0}A_i(2/\beta^4-4/(\ell\beta^3)) + 2\gamma/\xi_0 (-\ell(b_{1,i}+b_{2,i})/\beta^2 + b_{1,i}\beta)),
    \\& \mathcal{I}_i = -A_i(\ell^2/\beta),
    \\& \mathcal{J}_i = A_i(\ell^2/\beta^2 - 2\ell/\beta),
    \\& \mathcal{K}_i = 2A_i\ell/\beta^2,
    \\& \mathcal{L}_i = (-A_i\ell^2/\beta^2 + \ell(b_{1,i}+b_{2,i})),
    \\& \mathcal{M}_i = (2A_i\ell/\beta^2 + b_{1,i}),
    \\& \mathcal{N}_i = -(\avg{y_1F_i}+\avg{y_2F_i}) \ell^2,
    \\& \mathcal{O}_i = - (2\avg{y_1F_i}+\avg{y_2F_i}) \ell,
\end{align*}
\begin{align*}
    & \mathcal{P}_i = (\avg{y_1F_i} + \frac{1}{2} \avg{y_2F_i}),
    \\& \mathcal{R}_i = A_i(-\ell^2/\beta^2 + \ell/\beta),
    \\& \mathcal{U}_i = (A_i\ell/\beta^2 + b_{2,i}),
    \\& A_i = \gamma/\xi_0 \avg{y_1F_i},
    \\& b_{1,i} = 1/\alpha^2 (-(\alpha+\ell)\gamma/\xi_0(\avg{y_1F_i}+(1/2)\avg{y_2F_i}) - \ell^2\avg{xF_i}),
    \\& b_{2,i} = 1/\alpha^2 (\ell\gamma/\xi_0(\avg{y_1F_i}+(1/2)\avg{y_2F_i} - (k/\xi_0) \ell\avg{xF_i}),
\end{align*}
\end{fleqn}
and $
F_1 = ex, \: 
F_2 = \frac{2e}{\ell} \dot x, \: 
F_3 = \frac{e}{\ell^2} \ddot x, \:
F_4 = \frac{2e}{\ell} x, \: 
F_5 = \frac{e}{\ell^2} \dot x, \: 
F_6 = \frac{e}{\ell^2} x.
$

\section*{Acknowledgement}
L.M. is thankful for support through NSFC Grant No. 12271364 and GRF Grant No. 11302823.


\bibliography{references}
\bibliographystyle{unsrt}

\end{document}